\RequirePackage{rotating}
\documentclass[envcountsect,envcountsame]{svjour3}
\usepackage{algorithm}
\usepackage{graphicx}
\usepackage{mathtools}
\usepackage{mathrsfs}
\usepackage{fix-cm}
\usepackage{multirow}
\usepackage{amsmath}
\usepackage{amssymb}
\usepackage{latexsym}
\usepackage{dsfont}
\usepackage{xcolor}
\usepackage{cite}
\usepackage{dsfont}
\usepackage{enumitem}
\usepackage{todonotes}
\usepackage{hyperref}
\usepackage{booktabs}
\usepackage[group-separator={\,},group-minimum-digits=4]{siunitx}
\usepackage{accents}
\usepackage{adjustbox}
\usepackage{array}
\usepackage{subcaption}
\usepackage{caption}
\usepackage{multirow}

\usepackage{tabularray}
\UseTblrLibrary{siunitx}

\definecolor{noteD}{RGB}{211, 107, 145}

\usepackage[mathlines]{lineno}

\hypersetup{
	colorlinks=true,
	linkcolor=blue,
	citecolor=darkgreen,
	urlcolor=blue}
\definecolor{darkgreen}{RGB}{0,150,0}
\definecolor{C0}{RGB}{31,119,180}
\definecolor{C1}{RGB}{255, 127, 14}
\usepackage{algpseudocode}

\newtheorem{experiment}{Experiment}
\newtheorem{assumption}{Assumption}
\let\oldexperiment\experiment
\renewcommand{\experiment}{\oldexperiment\normalfont}

\def\tto{\;{\lower 1pt \hbox{$\rightarrow$}}\kern -10pt
	\hbox{\raise 2pt \hbox{$\rightarrow$}}\;}

\DeclareMathOperator{\dom}{dom}

\DeclareMathOperator{\prox}{prox}

\DeclareMathOperator*{\argmin}{argmin}

\newcommand{\Rex}{\overline{\mathbb{R}}}

\newcommand{\R}{\mathbb{R}}
\newcommand{\N}{\mathbb{N}}

\newcommand{\Ebb}{\mathbb{E}}

\newcommand{\Pb}{\mathbb{P}}

\newcommand{\toas}{\xrightarrow{\text{a.s.}}}

\newcommand\munderbar[1]{%
  \underaccent{\bar}{#1}}

	\definecolor{ceruleanblue}{rgb}{0.16, 0.32, 0.75}

\let\epsilon\varepsilon
\let\subseteq\subset

\newcommand{\varTh}{\mathrm{\Theta}}
\newcommand{\varLam}{\mathrm{\Lambda}}

\newcommand{\mU}{V}
\newcommand{\mV}{W}
\begin{document}
	
	\titlerunning{Randomized block proximal  method with  locally Lipschitz continuous gradient }
	\title{Randomized block proximal  method with  locally Lipschitz continuous  gradient
		\thanks{This work was partially supported by Conselleria d'Educació, Cultura i Universitats, Generalitat Valenciana  (grant number CIGE/2026/80). P. P\'erez-Aros  was partially supported by Centro de
			Modelamiento Matem\'{a}tico (CMM), ACE210010 and FB210005, BASAL funds for
			center of excellence and ANID-Chile grant: MATH-AMSUD 23-MATH-09 and
			MATH-AMSUD 23-MATH-17, ECOS-ANID ECOS320027, Fondecyt Regular 1220886,
			Fondecyt Regular 1240335, Fondecyt Regular 1240120 and Exploraci\'on  13220097.
			D. Torregrosa-Bel\'en was supported by Centro de Modelamiento Matemático (CMM) BASAL fund FB210005 for center of excellence from ANID-Chile and Fondecyt Postdoctorado 3250039, and by grants PID2022-136399NB-C21 and PID2025-168246NB-I00 funded by
ERDF/EU and by MICIU/AEI/10.13039/501100011033.}}
	\subtitle{}
	\author{\mbox{Pedro P\'erez-Aros} \and \mbox{David Torregrosa-Bel\'en} }
	
	\institute{Pedro P\'erez-Aros. ORCID: \href{https://orcid.org/0000-0002-8756-3011}{ 0000-0002-8756-3011} \at Departamento de Ingeniería Matemática and Centro de Modelamiento Matemático (CNRS IRL2807), Universidad de Chile, Beauchef 851, Santiago, Chile\\
		\email{pperez@dim.uchile.cl}
		\and Corresponding author: David Torregrosa-Bel\'en. ORCID: \href{https://orcid.org/0000-0003-2361-2037}{ 0000-0003-2361-2037}  \at Department of Mathematics, University of Alicante, Alicante, Spain\\ \email{david.torregrosa@ua.es}}
	
	\date{\today}
	
	\maketitle
	
	\begin{abstract}
	Block-coordinate algorithms are recognized to furnish efficient iterative schemes for addressing large-scale problems, especially when the computation of full derivatives entails substantial memory requirements and computational efforts. In this paper, we propose a randomized block proximal gradient algorithm for minimizing the sum of a smooth function and a separable proper lower semicontinuous function, both possibly nonconvex. In contrast to previous work, we only assume that the partial gradients of the smooth function are locally Lipschitz continuous with respect to their block of coordinates. At each iteration, the method adaptively selects a proximal stepsize to satisfy a sufficient decrease condition without requiring knowledge of the local Lipschitz moduli of the partial gradients of the smooth function. Hence, our work extends the applicability of randomized block-coordinate descent methods to more general problems where global Lipschitz gradient assumptions are not satisfied. In addition, our algorithm incorporates the possibility of conducting an additional boosted linesearch to enhance the performance of the algorithm. Our main result establishes subsequential convergence to a stationary point of the problem almost surely. 
	{Numerical experiments in nonnegative matrix factorization report that our method competes well with  recognized algorithms. Furthermore, a symmetric variant of nonnegative matrix factorization illustrates the advantage of the boosted linesearch in problems where the proximal stepsize cannot be tightly estimated in accordance to the local Lipschitz modulus.}
	\end{abstract}\vspace*{-0.05in}

	\vspace*{-0.05in}
	
	\vspace*{-0.2in}
	\section{Introduction}

Block-coordinate algorithms are traditionally among the most widely used approaches to tackle mathematical programs of diverse structures~(see, e.g.,~\cite{MR359795,MR273810,MR3444832,MR1099605,MR3347548,MR3719240} and the references therein). By iteratively updating different blocks of coordinates independently, these methods manage to reduce the complexity and memory requirements per iteration. This leads to recognized benefits for very-large-scale problems, especially when obtaining full derivatives  entails substantial computational efforts. Fields where the application of block-coordinate methods has proven fruitful are tomography, regularized least-squares, signal processing, nonnegative matrix factorization and protein structure prediction~\cite{491321,MR2421312,MR3158055,CDPR2003,Liu-2009-10238,06ae6e67-7fa9-32aa-97a2-a282f021e2c2}, while more can be consulted in~\cite{MR3347548}.

The block of coordinates to be updated at each iteration is typically determined using a cyclic, greedy, or randomized scheme.  Pathological examples are known to occur for the cyclic approach, as demonstrated by Powell~\cite{MR321541}, who showed that the cyclic coordinate gradient descent applied to a differentiable function can orbit around six points with nonzero gradient. On the other hand, greedy algorithms still require the computation of full derivatives and may present worse convergence rates than simple (1-block) schemes~\cite{MR2968857}.
The class of block-coordinate methods experienced a notable resurgence after the seminal work of Nesterov~\cite{MR2968857}, who established complexity results for the convergence of the randomized block-coordinate gradient descent in the convex case.  Nesterov's results were generalized to convex composite programs in~\cite{MR3179953}. Afterwards, the convex setting has been investigated by several works proposing parallel, distributed, inexact and accelerated variants (see, e.g.,~\cite{10.5555/3104482.3104523,MR3404687,MR3459207,MR3513271,MR3369495,MR3517098}).  {Furthermore, strategies for selecting random cycles in quadratic convex optimization have recently been studied in \cite{wright2020analyzing,gurbuzbalaban2020randomness}.}

This article investigates a randomized block-coordinate descent method adept at addressing  structured minimization problems in the form of
\begin{equation}\label{eq:P1}
    \min_{x\in\mathbb{R}^n} \varphi(x):=f(x)  + g(x), \tag{$\mathcal{P}$}
\end{equation}
where  the function $f\colon \R^n\to\R$ is smooth, $g\colon \R^n\to\Rex$ is    separable  and both are possibly nonconvex.
To the best of our knowledge, the first work dealing with the nonconvex framework is due to Patrascu and Necoara~\cite{MR3293507}. Specifically, they proposed  proximal gradient algorithms for problems in the form of~\eqref{eq:P1} when $g$ is assumed to be convex.
The same authors later considered the framework where $g$ is the $\ell_0$-norm and $f$ is convex~\cite{MR3365070}. The fully nonconvex setting is tackled in~\cite{MR3730541}.

The convergence analysis of the previously mentioned works (also in the convex setting) relies on a notion of block-wise global Lipschitz continuity of the gradient of the smooth $f$. Convergence guarantees depend on preselecting stepsizes satisfying 
 upper bounds  that involve  the inverse of the Lipschitz moduli of  block partial gradients. This compromises the performance of the method when the Lipschitz moduli cannot be estimated accurately. As an exception, Lu and Xiao~\cite{MR3730541} presented a proximal gradient method for~\eqref{eq:P1} that selects the proximal stepsize adaptively without assuming knowledge of the Lipschitz constants of the partial gradients. However, their analysis still requires the usual global Lipschitz gradients continuity assumption, which is needed to ensure that the sequence of stepsizes chosen by the adaptive procedure is bounded away from zero. The latter is fundamental in order to ensure almost sure subsequential convergence to stationary
 points of~\eqref{eq:P1}. 
 
 In this work, we go one step further by  dropping the global block-wise Lipschitz assumption and only requiring  the gradient of $f$ to be block-wise locally Lipschitz continuous, thus  extending  the applicability of our method to more general fields where global Lipschitz assumptions are violated, for instance, in nonnegative matrix factorization~\cite{MR2388467,LeeNMF,PAUCA200629}, quadratic inverse problems~\cite{MR3832977,MR3080197}, quadratic compressed sensing~\cite{MR3160324}, etc.
 Furthermore, from the theoretical point of view this demands developing a more sophisticated analysis  to ensure that an appropriate subsequence of stepsizes is bounded away from zero. 
 
We note that in the recent work~\cite{9414191}, only continuous differentiability of $f$ is assumed. This  framework is addressed through a \emph{randomized Bregman block proximal gradient} method. However, their analysis is based on the case where $g$ is convex and requires the selection of a family of auxiliary functions tailored to each problem that allows obtaining a sufficient decrease condition for $f$ block-wise. 
Therefore, the method we propose here is more general and at the same time allows for simpler  straightforward implementation.

In addition, as a complementary  contribution we permit the possibility of incorporating  in our numerical scheme an optional linesearch  in the spirit of  \emph{boosted  linesearches}. The term boosted linesearch was coined in the seminal works~\cite{MR4078808,MR3785672}, which proposed an accelerated variant for the classical  \emph{difference of convex functions algorithm}~\cite{MR874369} under some differentiability assumptions.  Afterwards, this  technique has been applied in different situations (see, e.g.,~\cite{MR4757561,MR4527568,tran2023boosteddcalgorithmclustering}). In the context of proximal  methods boosted linesearches have  been considered in the recent work~\cite{BDSA2025}, where they  were shown to enhance the performance of proximal (sub)gradient algorithms.  To the best of the authors' knowledge,   this is the first work that  introduces the use of boosted linesearch procedures in the  randomized coordinate framework.

{
\subsection{Related works in nonconvex proximal methods}

The research on proximal algorithms for the nonconvex setting has gained increasing attention after the seminal works of Attouch, Bolte and collaborators on the convergence properties of \emph{descent methods}~\cite{MR2421270,MR3010421}. Particularly, the works~\cite{MR3232623,MR2674728} analyze  proximal alternating methods which can be regarded as block proximal algorithms limited to two blocks. A landmark contribution is the \emph{iPiano algorithm} by Ochs et al.~\cite{MR3218822} which incorporates an inertial term in a proximal gradient iteration. A backtracking linesearch can be used to tune the inertial parameter when the Lipschitz modulus of the gradient is unknown. Another notable work is due to Bredies, Lorenz and Reiterer~\cite{MR3327417}, who investigated proximal methods  employing thresholding techniques to address fully nonconvex problems in \emph{Hilbert spaces}. More recently, Yang~\cite{MR4690351} proposed a proximal gradient scheme that combines extrapolation with a nonmonotone linesearch. In addition, some recent works have investigated the performance of proximal gradient methods beyond the restrictive global Lipschitz continuity assumption of the gradient. We highlight the approach by Bolte et al.~\cite{MR3832977}, that uses Bregman distances to derive an extended descent lemma inequality.  Meanwhile, other works exploit local Lipschitz gradient continuity by means of adaptive linesearches, see, e.g.,~\cite{8263933,MR4457955,MR4665669,MR4927930}.
}

\subsection{Main contributions}\label{maincontributions}
 We summarize the main contributions of the present work as follows:
\begin{itemize}

\item We propose a randomized block proximal gradient algorithm (see Algorithm~\ref{alg:1}) to address problems of the form~\eqref{eq:P1}, which does not require the partial gradients of the function \( f \) to be globally Lipschitz continuous. In addition, our method does not rely on prior knowledge of the (local) Lipschitz constants of the partial gradients of \( f \). At each iteration, a suitable proximal stepsize is selected through an adaptive procedure that ensures a sufficient decrease condition for~\( \varphi \). Moreover, a boosted linesearch can be incorporated to further enhance the algorithm’s performance.

\item In Theorem~\ref{t:2}, we demonstrate almost sure subsequential convergence of the sequence generated by our method to a stationary point of~\eqref{eq:P1}.   
The main difficulty is resolved in Proposition~\ref{prop:tau_bound}, where we prove that for bounded subsequences generated by Algorithm~\ref{alg:1} the infimum of the sequence of stepsizes is bounded away from zero.

\item We illustrate the performance of our method in  real-world applications in the fields of image compression {and image clustering.}
\end{itemize}

\subsection{Outline of the paper}
The remainder of the paper is structured as follows.  In Section~\ref{s:notation} we introduce the notation and some basic concepts on variational analysis and probability theory. Specifically, in Section~\ref{s:va} we  study  some  properties of the proximity operator that will be needed throughout the work. Our proposed algorithm, which we term as \emph{adaptive randomized block-proximal gradient}, is presented in Section~\ref{sect:alg}. Section~\ref{sect:convergence} is devoted to analyzing its convergence. In Section~\ref{sect:numerical}  we illustrate the performance of the method for some nonnegative matrix factorization and {symmetric nonnegative matrix factorization problems} that arise in different applications. Finally, conclusions and future directions of research are presented in Section~\ref{sect:conclusion}.

\section{Notation and preliminaries}\label{s:notation}

Throughout this paper, we work on the $n$-dimensional Euclidean space $\R^n$. The notation $\langle \cdot, \cdot \rangle$ and $\| \cdot\|$ stand for the standard inner product and the $\ell_2$-norm in $\R^n$, respectively. The nonnegative reals are denoted as $\R_+$, while the extended-real-valued line is the set $\Rex:=\R\cup\{-\infty,+\infty\}$ and by convention we set $1/0=+\infty$. We use $\mathbb{S}^n$ for the space of $n\times n$ symmetric matrices while $\N:=\{0,1,2,\ldots\}$ is the set of nonnegative integers.  We use $A\succeq 0$ to indicate that a  symmetric matrix $A$ is positive semidefinite. In particular, if $A$ is diagonal this is equivalent to all its diagonal entries  belonging to $\R_{+}$. Given $A,B\in\mathbb{S}^n$, the expression $A \succeq B$ is equivalent to $A-B\succeq 0$.   We use $\to$ to denote convergence of a sequence and $\mathbb{B}_{\epsilon}(\bar{x}) := \{ x \, : \, \| x- \bar{x}\| < \epsilon\}$ for the open ball of radius $\epsilon >0$ centered at $\bar{x}$. 

In the following, we will split a vector  $x\in\R^n$ into $N\in{\mathbb{N}{\setminus}\{0\}}$ blocks of coordinates of dimension $n_i$, for $i\in{\{1,\ldots,N\}}$, i.e., $\sum_{i=1}^N n_i=n$. In order to access such blocks  it will be handy to consider a decomposition of the $n\times n$ identity matrix $I_n=[U_1, U_2,\ldots,U_N]$, where the number of columns of $U_i$ is $n_i$. Thus, the $i$-th  block of coordinates of a vector $x\in\R^n$ is obtained as $x_i:=U_i^Tx\in\R^{n_i}$.
Let $f\colon\R^n\to\R$ be a differentiable function. Similarly, we denote the $i$-th block partial gradient of $f$ at a point $x$ as $\nabla_i f(x) :=      U_i^T \nabla f(x)$, where $\nabla f$ is the gradient of $f$. A function is said to be \emph{smooth} if it is  differentiable with continuous gradient.

Given some constant $L\geq 0$, a vector-valued function $F\colon C\subseteq \R^n\to\R^m$ is said to be \emph{$L$-Lipschitz continuous} on~$C$ if
\[
\|F(x)-F(y)\| \leq L \|x-y\|, \quad \text{for all } x,y\in C,
\]
and \emph{locally Lipschitz continuous} around $\bar{x} \in C$ if it is $L$-Lipschitz continuous in some neighborhood of $\bar{x}$ for some $L\geq0$. In particular, $F$ will be said to be locally Lipschitz continuous on $C$ if it is locally Lipschitz continuous around every point in~$C$. 

It is well-known that if $f\colon  \R^n\to\R$  has $L$-Lipschitz continuous gradient on a convex set $C$  it verifies the so-called \emph{descent lemma} on  $C$, namely,
\[
f(y) \leq f(x) + \langle \nabla f(x), y-x  \rangle + \frac{L}{2} \|y-x\|^2, \quad \text{for all } x,y\in C.
\]

 {Now, we introduce some foundational concepts from probability theory and random sequences that are essential for analyzing the numerical methods with stochastic components considered in this paper. For further details, the reader is referred to the classical monographs~\cite{MR2267655,MR2722836,MR270403}.

Throughout the paper, we denote a (sample) probability space by $(\Omega, \mathcal{A}, \mathbb{P})$. For random variables defined on $\Omega$, we omit the explicit dependence on the sample point $\omega \in \Omega$ to simplify notation. The $\sigma$-algebra generated by a collection $\Phi$ of random variables is denoted by $\sigma(\Phi)$.

The expectation of a random variable $\phi \colon  \Omega \to \R$ is denoted by $\mathbb{E}[\phi]$. Given a sub-$\sigma$-algebra $\mathcal{F} \subseteq \mathcal{A}$, the conditional expectation of $\phi$ with respect to $\mathcal{F}$ is written as $\mathbb{E}[\phi \mid \mathcal{F}]$.

A sequence of random variables $(x^k)_{k \in \mathbb{N}}$ is said to converge to a random variable $x$ \emph{almost surely} (abbreviated \emph{a.s.}) if
\[
\mathbb{P}\left( \lim_{k \to \infty} x^k = x \right) = 1.
\]
In this case, we write $x^k \xrightarrow{\text{a.s.}} x$.}

\subsection{Preliminary notions on variational analysis}\label{s:va}

We now briefly introduce some concepts on variational analysis that will be required throughout the manuscript. Further details can be consulted in the classical references~\cite{MR2191744,MR3823783,MR1491362}.

Let $g\colon \R^n\to\Rex$ be an extended-real-valued function, and define its \emph{domain} as $\dom{g}:=\{x\in\R^n\,:\, g(x) < +\infty\}$. We say that $g$ is proper if it does not attain the value $-\infty$ and $\dom{g}\neq\emptyset$. The function $g$ is \emph{lower semicontinuous (lsc)} at some point $\bar{x}\in\R^n$ if $\liminf_{x\to\bar{x}} g(x) \geq g(\bar{x})$.

Given $\bar{x}\in  g^{-1}(\R)$, the \emph{regular subdifferential} of the function $g$ at $\bar{x}\in\mathbb{R}^n$ is the closed and convex set of \emph{(regular) subgradients}
 \begin{equation*}
 \partial g(\bar{x}) := \bigl\{ v\in\mathbb{R}^n : g(x) \geq g(\bar{x}) + \langle v, x-\bar{x}\rangle + o(\|x-\bar{x}\|) \bigr\}.
\end{equation*}
We use the convention $ {\partial} g(\bar{x}):=\emptyset$ if $|g(\bar x) |= +\infty$. 
Recall that for convex functions the  regular subdifferential is equal to the standard convex subdifferential (see, e.g.,~\cite[Definition~16.1]{bauschke2017}). {By~\cite[Proposition~1.30(i)]{MR3823783}} a necessary condition for a point $x^*\in  \varphi^{-1}(\R)$ to be a local minimum of $\varphi$ is that
\begin{equation}\label{eq:stationary_point}
0 \in \partial \varphi(x^*). 
\end{equation}
A point $x^*$ satisfying~\eqref{eq:stationary_point} is said to be a \emph{stationary point} of~\eqref{eq:P1}. In particular, for a function $\varphi$ defined as in~\eqref{eq:P1} the calculus rules of the subdifferential (see, e.g.,~\cite[Proposition~1.30(ii)]{MR3823783}) yield
\begin{equation*} 
0 \in \partial \varphi(x^*) \Leftrightarrow 0 \in  \nabla f(x^*) + \partial g(x^*).
\end{equation*}


We define the \emph{proximal point mapping} (also called \emph{proximity operator}) of a proper lsc function $g\colon \R^n\to\Rex$ as the multifunction  $\prox_g\colon \R^n\tto\R^n$ defined as
\[
\prox_{g}(x):=\argmin_{u\in\R^n} \left\{ g(u) + \frac{1}{2} \|x-u\|^2 \right\}.
\]
A function $g\colon \R^n\to\Rex$ is said to be \emph{prox-bounded} if there exists some $\tau >0 $ such that $g(\cdot)  + \frac{1}{2\tau}\|\cdot\|^2$ is bounded from below { on $\R^n$}. The supremum of the set of all such constants is called the \emph{prox-boundedness threshold} of $g$ and is denoted by $\tau^{g}$. 
For a proper, lsc and prox-bounded function, the proximal mapping $\prox_{\tau g}$ has full domain for any  $\tau\in{]0, \tau^{g}[}$, but it may not be single-valued (see, e.g., \cite[Theorem 1.25]{MR1491362}). If $g$ is proper, lsc and convex then $\tau^{g} = +\infty$ and $\prox_{\tau g}$ is single-valued for all $\tau>0$. The following lemmas gather properties of prox-bounded functions that may be needed in the sequel.  The first lemma establishes the stability of prox-boundedness and the proximal operator under linear perturbations of the function.
 Given a function $g\colon \R^n\to\Rex$ and $v\in \mathbb{R}^n$, we denote 
\begin{align}\label{def_varphiv}
	g^{v} := g + \langle \cdot, v \rangle.
\end{align}
\begin{lemma}\label{lemma:proxv}
Let $g\colon \R^n\to\Rex$ be a proper lsc prox-bounded function with prox-boundedness threshold $\tau^{g}>0$. Given $v\in\R^n$,   the function $g^{v}$  is   prox-bounded with $\tau^{g^{v}}=\tau^{g}$. In addition, $\prox_{\tau g^{v}}(\cdot) = \prox_{\tau g}(\cdot - \tau v)$, for all $\tau \in{]0,\tau^{g}[}$.
\end{lemma}
\begin{proof}
{ 
We start by showing that $g^{v}(\cdot) + \frac{1}{2\tau}\|\cdot\|^2$ is bounded from below on $\R^n$ for all $\tau\in{]0,\tau^g[}$. Indeed, let  $\tau\in{]0,\tau^g[}$ and $\varepsilon>0$. \emph{Cauchy--Schwarz's} and \emph{Young's inequalities} yield
\[
\begin{aligned}
g^v(x) + \frac{1}{2\tau} \|x\|^2  & = g(x) + \langle x, v\rangle + \frac{1}{2\tau}\|x\|^2  \\
& \geq g(x) + \frac{1}{2} \left(\frac{1}{\tau} - \frac{1}{\epsilon} \right) \|x\|^2 - \frac{\varepsilon}{2}\|v\|^2,
\end{aligned}
\]
for all $x\in\R^n$. By setting $\varepsilon > (1/\tau - 1/\tau^{g})^{-1}$, we get that $(1/\tau - 1/\varepsilon) > 1/\tau^{g}$ and the prox-boundedness of $g$ confirms the claim. Hence, we have proven that $g^v$ is prox-bounded with prox-boundedness threshold $\tau^{g^v} \geq \tau^g$. To show that $\tau^{g^{v}}=\tau^{g}$, observe that  $g=g^{v} + \langle -v,x\rangle$. Thus, the above implies that $g(\cdot) + \frac{1}{2\tau}\|\cdot\|^2 $ is bounded from below on $\R^n$ for all $\tau\in{]0,\tau^{g^v}[}$, hence $\tau^{g_v}\leq \tau^g$.} 

Furthermore, let $\tau\in{]0,\tau^{g}[}$ and fix $x\in\R^n$. For all $u\in\R^n$, we have 
\[
 g(u) + \frac{1}{2\tau}\|u-(x-\tau v)\|^2 = g(u) + \frac{1}{2\tau}\|u-x\|^2 + \langle u,v\rangle  - \langle x,v\rangle+ \frac{\tau}{2} \|v\|^2. 
 \]
Therefore, 
\[
\begin{aligned}
\argmin_{u\in\R^n}  \left\{ g(u) + \frac{1}{2\tau}\|u-(x-\tau v)\|^2 \right\} 
& =  \argmin_{u\in\R^n}  \left\{ g(u) + \frac{1}{2\tau}\|u-x\|^2 + \langle u,v\rangle \right\},
\end{aligned}
\]
which proves the relationship between the proximity operators of $g$ and $g^{v}$.
\end{proof}

The second lemma establishes a local property of outer semicontinuity of the proximal operator with respect to the parameter $  \tau > 0$. More precisely, we have the following result
\begin{lemma}\label{lemma:proxclosed}
Let $g\colon \R^n\to\Rex$ be a proper lsc prox-bounded function with prox-boundedness threshold $\tau^{g}>0$. Let  $\bar{x}\in \dom g$ and $\epsilon>0$. Then there exists a constant $\tilde{\tau} \in{]0,\tau^{g}[}$ such that $x\in\mathbb{B}_{\epsilon}(\bar{x})$ for all $x\in\prox_{\tau g}(\bar{x})$ and $\tau < \tilde{\tau}$.
\end{lemma}

\begin{proof}
By way of contradiction, let us assume there exists a sequence  $(\tau_k)_{k\in\N} \subset {]0,\tau^{g}[}$ with $\tau_k\to 0^+$ and such that there exists
\[
x^k \in \prox_{\tau_k g}(\bar{x}) \quad \text{with} \quad \|x^k-\bar{x}\| \geq \epsilon, \text{ for all } k\geq 0.
\]
 By the definition of the proximity operator, this yields
\begin{equation}\label{eq:prox-lemma-1}
g(x^k) + \frac{1}{2\tau_k} \|x^k-\bar{x}\|^2 \leq g(\bar{x}), \quad \forall k\geq 0.
\end{equation}
Now, fix  $\hat{\tau}\in{]0,\tau^{g}[}$ { and let $\delta> (1/\hat{\tau}-1/\tau^g)^{-1}$. By resorting to \emph{Cauchy--Schwarz's} and \emph{Young's inequalities} we get
\begin{equation}\label{eq:p_CSY}
\begin{aligned}
g(x^k) + \frac{1}{2 \hat{\tau}} \|x^k - \bar{x}\|^2 & = g(x^k) + \frac{1}{2\hat{\tau}}\|x^k\|^2 - \frac{1}{\hat{\tau}} \langle x^k, \bar{x}\rangle + \frac{1}{2\hat{\tau}} \|\bar{x}\|^2 \\
& \geq g(x^k) + \frac{1}{2} \left(\frac{1}{\hat{\tau}} - \frac{1}{\delta} \right) \|x^k\|^2   + \frac{1}{2} \left(\frac{1}{\hat{\tau}}-  \frac{\delta}{\hat{\tau}^2} \right) \|\bar{x}\|^2,
\end{aligned}
\end{equation}
for all $k\geq 0$.} Thus, the prox-boundedness of $g$ implies that there exists a constant $\beta\in\R$ such that 
\begin{equation}\label{eq:prox-lemma-2}
g(x^k) \geq \beta - \frac{1}{2\hat{\tau}} \|x^k -\bar{x}\|^2, \quad \forall k\geq 0.
\end{equation}
Putting~\eqref{eq:prox-lemma-1}-\eqref{eq:prox-lemma-2} together and dividing by $\|x^k-\bar{x}\|^2$   we have
\[
\frac{1}{2\tau_k} - \frac{1}{2\hat{\tau}} \leq \frac{g(\bar{x}) - \beta}{\|x^k-\bar{x}\|^2} \leq \frac{g(\bar{x})-\beta}{\epsilon^2}, \quad \forall k\geq 0.
\]
Since the left-hand side of the above expression tends to $+\infty$ as $\tau_k\to0^+$ we get a contradiction for sufficiently large $k$. This proves the result.
\end{proof}

The final lemma establishes the closedness of the graph of the proximal mapping with respect to all of its parameters.

\begin{lemma}\label{lemma:proxcoledk}
Let $g\colon \R^n\to\Rex$ be a proper lsc  prox-bounded function with prox-boundedness threshold $\tau^{g}>0$. Assume $x^k \to \bar{x}\in\dom{g}$, $\tau_{k}\to0^+$, {with $\tau_k< \tau^{g}$ for all $k$,} and let $w^k \in \prox_{\tau_{k} g} (x^k)$. Then $w^k \to \bar{x}$.
\end{lemma}

\begin{proof}
By the definition of proximity operator we have
\begin{equation}\label{eq:lm2-0}
g(w^k)+ \frac{1}{2\tau_k}\|w^k-x^k\|^2 \leq g(\bar{x}) + \frac{1}{2\tau_k} \|\bar{x}- x^k\|^2.
\end{equation}
{ Fix $\hat{\tau}\in{]0,\tau^g[}$. By proceeding similarly to~\eqref{eq:p_CSY}, but replacing $x^k$ by $w ^k$ and $\bar{x}$ by $x^k$, we get 
\[
g(w^k) + \frac{1}{2 \hat{\tau}} \|x^k - w^k \|^2 \geq g(w^k) + \frac{1}{2} \left(\frac{1}{\hat{\tau}} - \frac{1}{\delta} \right) \|w^k\|^2   + \frac{1}{2} \left(\frac{1}{\hat{\tau}}-  \frac{\delta}{\hat{\tau}^2} \right) \|x^k\|^2, \quad \forall k\geq 0,
\]
where again $\delta > (1/\hat{\tau}-1/\tau^g)^{-1}$.
Hence, by using the prox-boundedness of $g$ and the fact that $(x^k)_{k\in\N}$ is bounded,  we can ensure that there exists    $\beta\in\R$ such that
\begin{equation}\label{eq:lm2-1}
g(w^k) \geq \beta - \frac{1}{2\hat{\tau}} \|w^k-x^k\|^2, \quad \forall k\geq 0.
\end{equation}}
Finally, combining \eqref{eq:lm2-0} and \eqref{eq:lm2-1}, we get that for large enough $k \in\N$ it holds
\begin{equation*}
  \|w^k-x^k\|^2 \leq   \left( \frac{1}{2} -  \frac{\tau_k}{2\hat{\tau}}\right)^{-1} \left( \tau_k \bigl(g(\bar{x}) - \beta\bigr) + \frac{1}{2} \|\bar{x}- x^k\|^2 \right),
\end{equation*}
which shows that  $\|w^k-x^k\| \to 0$, and consequently $w^k \to \bar x$.  
\end{proof}

\section{Adaptive randomized block proximal gradient algorithm}\label{sect:alg}

In this work, we are interested in the structured minimization problem \eqref{eq:P1} when  the functions $f\colon \R^n\to\R$  and $g\colon \R^n\to\Rex$ satisfy the conditions gathered in the following assumption.
 
\begin{assumption}\label{assumption:1} {Let  all vectors $x\in\R^n$ be decomposed into $N\in{\mathbb{N}{\setminus}\{0\}}$ blocks  of coordinates $x_i\in\mathbb{R}^{n_i}$ of dimension $n_i$, for $i\in{\{1,\ldots,N\}}$, with $\sum_{i=1}^N n_i=n$.} We assume that:
\begin{enumerate}[label=(\roman*)] 
\item\label{it:as-1} $f\colon \R^n\to\R$ is smooth and has block-wise  locally Lipschitz continuous gradient on $\dom g$, i.e., for  every point $\bar{x}\in \dom g$ there exist {convex} neighborhoods $V_i$ and constants $L_i(\bar{x})>0$ such that, for all $i\in{\{1,\ldots,N\}}$, it holds
\begin{equation}\label{eq:descent_lemma_0}
\| \nabla_i f(x+U_i s_i) - \nabla_i f(x) \| \leq L_i(\bar{x}) \| s_i \|, \quad \forall x,x + U_is_i \in V_i,
\end{equation}
 where the constants $L_i(\bar{x})$ depend  both  on the point $\bar{x}$ and the block-coordinate~$i$;
\item $g\colon \R^n\to\Rex$ is block-separable as
\[
 g(x)=\sum_{i=1}^N g_i(x_i), \quad \forall x\in \R^n,
\]
where the functions $g_i:\R^{n_i}\to\Rex$ are proper, lsc and prox-bounded for all $i\in\{1,\ldots,N\}$ {and are only evaluated in the $i$-th block $x_i\in\R^{n_i}$ of $x\in\R^n$};
 
\item $\inf_{x\in\R^n} \varphi(x) >-\infty$.

\end{enumerate}
\end{assumption}

\begin{remark}
Observe that Assumption~\ref{assumption:1}\ref{it:as-1} encompasses the standard assumption for block-coordinate methods of having a block-wise (globally) Lipschitz continuous gradient (see, e.g.,~\cite{MR2968857,MR3293507}). The latter is directly recovered from~\eqref{eq:descent_lemma_0} by taking for all $i\in\{1,\ldots,N\}$, $V_i=\mathbb{R}^n$ and $L_i(\bar{x}):=L_i$ for all $\bar{x}$. Nonetheless, having a block-wise locally Lipschitz continuous gradient is strictly more general, as demonstrated by the simple two-variable function $f(x_1,x_2)=x_1^2x_2^2$. This distinction also arises in the nonnegative matrix factorization model considered in our numerical tests below.
\end{remark}

By the descent lemma,~\eqref{eq:descent_lemma_0} yields,  for all $i\in{\{1,\ldots,N\}}$, 
\begin{equation}\label{eq:descent_lemma}
f(x+U_i s_i) \leq f(x) + \langle \nabla_i f(x), s_i\rangle +\frac{L_i(\bar{x})}{2} \|s_i\|^2, \quad \forall x,x + U_is_i \in V_i,
\end{equation}
provided that $V_i$ is convex. This inequality can be regarded as a block-wise local version of the descent lemma and is the core tool used in our analysis.

In Algorithm~\ref{alg:1} we present our proposed scheme for addressing~\eqref{eq:P1} under Assumption~\ref{assumption:1}.

 \begin{algorithm}[h!]
	\caption{Adaptive Randomized  Block Proximal Gradient  (ARBPG) algorithm}\label{alg:1}
	\begin{algorithmic}[1]
		\Require{Choose $x^0\in\dom{\varphi}$;  $0<\munderbar{\tau}\leq\bar{\tau}$; $a,\alpha>0$; $\beta\in{]0,1[}$. Set $k:=0$.}
		\State{Choose randomly a block of coordinates $i_k$ with probability {$\Pb(i{=}i_k){=}p_{i}>0$}.}
		\State{Set $\bar{\tau}_k\in{[\munderbar{\tau},\min\{\bar{\tau},\tau^{g_{i_k}}\}[}$. { Take $t=0$.} }
		\State{{Set $\tau_k := \beta^t\bar{\tau}_k$ and} compute a direction $d^k\in\mathbb{R}^n$  {with $d^k_{i}=0$, for all $i\neq i_k$, and $d^k_{i_k}$ given by the block-proximal step}
			\begin{equation}\label{eq:prox}
				\begin{aligned}
				d^k_{i_k} \in    \prox_{ \tau_{k} g_{i_k}}\left(   x^k_{i_k} - \tau_k \nabla_{i_k } f(x^k)  \right) - x^k_{i_k}.
				\end{aligned}
			\end{equation}}
		\State{\textbf{if} $d_{i_k}^k = 0$ \textbf{then} update   $x^{k+1}:=x^k$, $k:=k+1$ and go to Step~1.}
		\State{\textbf{if} $d^k$ verifies
		\begin{align}\label{eq:NMLS}
		\varphi(x^k+d^k)  \leq \varphi(x^k) - a\|d^k\|^2,
		\end{align}    
		\hspace{\algorithmicindent}\textbf{then} set $\hat{x}^{k}:=x^k+d^k$ and go to Step~6 		\textbf{else} {update $t\leftarrow t+1$} and go back to Step~3.}
		\State{Pick any $\lambda_k\geq 0$ such that}
		\begin{align}\label{eq:BLS}
		\varphi(\hat{x}^k+\lambda_k d^k)  \leq  \varphi(\hat{x}^k) - \alpha \lambda_k^2 \|d^k\|^2.
		\end{align}    
\State{Update $x^{k+1}=  \hat{x}^k + \lambda_k d^k$. Set $k:=k+1$ and go back to Step~1.}
	\end{algorithmic}
\end{algorithm}

{Every iteration of Algorithm~\ref{alg:1} starts by randomly selecting a block index according to some prescribed probability distribution  $\Pb(i{=}i_k){=}p_{i}>0$ (not necessarily uniform)  fixed for every step of the algorithm}. In Step~2   an arbitrary stepsize $\tau_k=\bar{\tau}_k$  is picked in order to perform the block proximal gradient step described by~\eqref{eq:prox}. Since this proximal stepsize may be arbitrarily taken with no information about the local Lipschitz constant of $f$, the algorithm includes  a linesearch  to ensure that the proximal step furnishes a sufficient decrease of the objective function,  which is fundamental for carrying out the convergence analysis of the method. Hence, Step~5 defines   an inner loop whose goal is to reduce the stepsize $\tau_k$ by a backtracking procedure until the sufficient decrease condition given by~\eqref{eq:NMLS} is satisfied.   {Note that if at some step a direction $d^k=0$ is taken, then by Step~4 the next iteration is updated as $x^{k+1} :=x^k$ and the sufficient decrease condition~\eqref{eq:NMLS} is trivially satisfied as an equality.}

Furthermore, we allow the possibility of incorporating an optional  linesearch  in Step 6 that follows the spirit of the \emph{boosted linesearches} in the works~\cite{BDSA2025,MR4078808,MR3785672}. The underlying idea of this linesearch is to examine whether by advancing further in the direction determined by the proximal step a larger decrease can be achieved. Note that the computational burden of  this step can be omitted by simply taking $\lambda_k=0$. In such a case, we update the iteration as $x^{k+1}:= \hat{x}^k$. Nonetheless, different works report that the implementation of an appropriate strategy for selecting $\lambda_k>0$ can lead to numerical improvements. We provide further details on this in what follows.

Observe that~\eqref{eq:BLS}  can be regarded as a linesearch for the function $\varphi$ in the interval $\{ \hat{x}^k + \lambda d^k \, : \, \lambda \geq 0 \}$. Hence, an essential fact to consider when addressing this step is whether the direction $d^k$ defines a descent direction for~$\varphi$ at the point~$\hat{x}^k$. In order to investigate this scenario, it is relevant to note that when only one block of variables is considered and the gradient of $f$ is globally Lipschitz continuous with known Lipschitz constant (and hence the adaptive procedure in Step~5 is not needed), Algorithm~\ref{alg:1} recovers a particular case of the \emph{Boosted Double-proximal Subgradient Algorithm} in~\cite{BDSA2025}. In view of this, it is easy to extend~\cite[Proposition~3.5]{BDSA2025} to prove that if $g$ is differentiable at~$\hat{x}^k$ then for any $\alpha>0$ there exists some $\gamma_k>0$ such that
\[
\varphi(\hat{x}^k + \lambda d^k)  \leq \varphi(\hat{x}^k) - \alpha \lambda^2 \|d^k\|^2, \quad \text{for all } \lambda\in{[0,\gamma_k]}.
\]
Thus, if $g$ is differentiable  one can conduct a backtracking-type linesearch in Step~6 of Algorithm~\ref{alg:1} to obtain a suitable $\lambda_k$ that satisfies~\eqref{eq:BLS}.

Nonetheless, the authors in~\cite{BDSA2025} also explored how to efficiently apply the boosted procedure in the absence of differentiability. More precisely, when $g$ is not differentiable they proposed to carry out only a fixed number of attempts for the linesearch. If the decrease condition is not satisfied after this number of attempts then the linesearch is stopped by setting $\lambda_k=0$. This strategy resulted very fruitful  for the non-differentiable setting and significantly outperformed the version of the  method without linesearch, as illustrated in the numerical experiments in~\cite{BDSA2025}. In Section~\ref{sect:numerical} below, we devise a new procedure  for   implementing boosted linesearches tailored for non-differentiable functions  and that does not require prespecifying a fixed  number of attempts in the linesearch.  

Before delving into  the convergence analysis of Algorithm~\ref{alg:1}, we present the following lemma, which ensures that the inner loop determined by Step~5 of Algorithm~\ref{alg:1} must terminate finitely.
\begin{proposition}\label{prop:ls_wd}
Suppose Assumption~\ref{assumption:1} holds. Let $a>0$ and $k \geq 0$  fixed. Then there exists $\zeta_k>0$ such that  for any $\tau_k\in{]0,\zeta_k[}$ and $s\in\R^n$ defined as
\begin{equation}\label{eq:ls_wd}
s_{i_k} \in    \prox_{ \tau_{k} g_{i_k}}\left(   x^k_{i_k} - \tau_k \nabla_{i_k } f(x^k)  \right) - x^k_{i_k},
\end{equation}
and $s_j=0$, for $j\neq i_k$,  the decrease condition
\[
\varphi(x^k + s) \leq \varphi(x^k) - a \|s\|^2
\]
holds.       
\end{proposition}

\begin{proof}
{For simplicity of notation, and since $k$ and $i_k$ are fixed, we denote $v := \nabla_{i_k} f(x^k)$. By Lemma~\ref{lemma:proxv}, we deduce from~\eqref{eq:ls_wd} that 
\begin{equation}\label{eq:proxxkdk}
x_{i_k}^k + s_{i_k} \in \prox_{\tau_k g_{i_k}^{v}} (x_{i_k}^k), 
\end{equation}
where we recall   that $g_{i_k}^{v}$ is defined in \eqref{def_varphiv}. Moreover, thanks to Lemma \ref{lemma:proxv},  the function $g_{i_k}^{v}$ has prox-boundedness threshold $\tau^{g_{i_k}}$. Let $V_{i_k}$ be a neighborhood of $x^k$  such that~\eqref{eq:descent_lemma}  holds with $\bar{x}:= x^k$ and $i:=i_k$ for some local Lipschitz constant  $L_{i_k}(x^k)>0$. By~\eqref{eq:proxxkdk} and Lemma~\ref{lemma:proxclosed} applied  to $g_{i_k}^{v}$, there exists $\tilde{\tau}\in{]0,\tau^{g_{i_k}}[}$ such that  $x_{i_k}^k + s_{i_k} \in V_{i_k}$ for all $\tau_k < \tilde{\tau}$.

Set $\zeta_k:=\min{\{\tilde{\tau},(L_{i_k}(x^k)+2a)^{-1}\}}$. Hence, for any $\tau_k < \zeta_k$, applying the descent lemma we get
\begin{equation}\label{eq:sat_ls}
\begin{aligned}
\varphi(x^k+s)& = f(x^k+U_{i_k} s_{i_k}) + g(x^k+U_{i_k} s_{i_k}) \\
&\leq f(x^k) + g(x^k+U_{i_k} s_{i_k}) +  \langle \nabla_{i_k} f (x^k), s_{i_k} \rangle + \frac{L_{i_k}(x^k)}{2} \|s_{i_k}\|^2 \\
&= \varphi (x^k)  + g_{i_k}(x_{i_k}^k+s_{i_k}) +  \langle \nabla_{i_k} f(x^k), s_{i_k} \rangle + \frac{1}{2\tau_k}\|s_{i_k}\|^2 - g_{i_k}(x_{i_k}^k) \\ & \quad   + \left(\frac{L_{i_k}(x^k)}{2} - \frac{1}{2\tau_k} \right)\|s_{i_k}\|^2 \\
&\leq  \varphi (x^k) - a\|s_{i_k}\|^2,
\end{aligned}
\end{equation}
where the  last inequality is due to \eqref{eq:proxxkdk} and $\tau_k \leq 1/(L_{i_k}(x^k)+2a)$. This completes the proof.}
\end{proof}

{
\begin{remark}[On the choice of the trial stepsize]\label{remark:trial_step}
Equation~\eqref{eq:sat_ls} implies that  a candidate stepsize $\tau_k$ will satisfy the sufficient decrease condition~\eqref{eq:NMLS} whenever 
\begin{equation}\label{eq:sat_ls2}
   \tau_k \leq  \frac{1}{L_{i_k}(x^k)+2a}.
\end{equation}
Consequently, \eqref{eq:sat_ls2} highlights that the trial stepsize parameter $\bar{\tau}_k$ crucially influences algorithmic performance. If $\bar{\tau}_k$ is significantly larger than $1/L_{i_k}(x^k)$, several Armijo linesearch iterations may be required to gurantee sufficient decrease, thereby increasing the computational burden. Conversely, if $\bar{\tau}_k \ll 1/L_{i_k}(x^k)$, the linesearch condition is easily satisfied, but the resulting stepsize may be excessively small, leading to slow convergence. 
Ideally, $\bar{\tau}_k$ should approximate the reciprocal of the local Lipschitz modulus. However, computing the exact Lipschitz constant is often computationally expensive. In practice, it is more manageable to choose $\bar{\tau}_k$ based on an estimate or upper bound of $L_{i_k}(x^k)$. In this scenario, a boosted linesearch in the spirit of~\eqref{eq:BLS} becomes particularly fruitful, as it allows for larger decreases in the objective function without the extra cost of evaluating new derivatives. We investigate the numerical benefits of this boosted approach in Section~\ref{sect:num3}.
\end{remark}
}
\section{Convergence analysis}\label{sect:convergence}

In what follows, we study the convergence properties of the random sequences generated by Algorithm~\ref{alg:1}.
Observe that at iteration $k$, Algorithm~\ref{alg:1}  generates   random vectors, namely $\bigl(x^k, d^{k-1},\lambda_{k-1}\bigr)$, defined on the probability space $(\Omega,\mathcal{A}, \mathbb{P})$. {Note that Algorithm~\ref{alg:1} requires the selection of the initial trial stepsize in Step~3 and one of possibly multiple solutions of the proximal subproblem~\eqref{eq:prox}. Hence, to ensure the measurability of the random vectors $\bigl(x^k, d^{k-1},\lambda_{k-1}\bigr)$ one can assume the existence of a deterministic selector for these operations. Under this selection rule, the sequence $\bigl(x^k, d^{k-1},\lambda_{k-1}\bigr)$ becomes an adapted process with respect to the filtration defined by the random indices, which is a fundamental requirement for our main results. Instead of imposing restrictions on the implementation of Algorithm~\ref{alg:1} (such as  uniqueness of the solutions of the proximal subproblems) in what follows we shall assume that the random sequences generated by Algorithm~\ref{alg:1} are measurable.}


{Our first convergence result for Algorithm~\ref{alg:1} shows that the sequence of function values $\bigl(\varphi(x^k)\bigr)_{k\in\mathbb{N}}$ converges almost surely to the random variable given by its infimum. Moreover, we prove that the sequence of squared distances between consecutive iterates is almost surely summable. To establish this latter fact, we note that when the optional linesearch in Step~6 of Algorithm~\ref{alg:1} is used, it is not generally true that $\|d^k\| = \|x^{k+1} - x^k\|$. }The presence of the additional random sequence $(\lambda_k)_{k\in\N}$ leads to the alternative relation $x^{k+1}-x^k = (1+\lambda_k)  d^k$, for $k\geq 0$. Nevertheless, in the following result we show that it is possible to reformulate the sufficient decrease  condition resulting from Algorithm~\ref{alg:1} in terms of this   sequence (see \eqref{eq:suff_dec_0} below). 

\begin{theorem}\label{t:1}
Suppose that Assumption~\ref{assumption:1} holds.  Given $x^0\in\dom \varphi$, consider the sequence $(x^k)_{k\in\N}$ generated by Algorithm~\ref{alg:1}. Then, for all $k\geq 0$,
\begin{equation}\label{eq:dec1}
\varphi(x^{k+1}) \leq \varphi(x^k) - \frac{\eta}{2} \|x^{k+1}-x^k\|^2,
\end{equation}
where    $\eta:=\min\{\alpha,a\}$.  In addition, the following assertions hold:
\begin{enumerate}[label=(\roman*)] 
\item\label{it:t1-1} The sequence $\bigl(\varphi(x^{k})\bigr)_{k\in\N}$ is monotonically nonincreasing and converges almost surely to the real random variable $\varphi^*_{\infty}:=\inf_{k\in\N} \varphi(x^{k})$. Furthermore, $\lim_{k\to\infty} \Ebb \bigl[\varphi(x^{k})\bigr] =  \Ebb \bigl[\varphi^*_{\infty}\bigr]$.
\item\label{it:t1-2} It holds that
\begin{equation}\label{eq:sum_deltak}
\mathbb{E}\left[ \sum_{k=0}^{\infty} \|x^{k+1}-x^k\|^2 \right]< \infty.
\end{equation}
Particularly, $\|x^{k+1}-x^k\|\toas 0$ and $\Ebb[\|x^{k+1}-x^k\|]\to 0$ as $k\to\infty$.
\end{enumerate}
\end{theorem}
\begin{proof}
Let $\eta=\min{\{\alpha,a\}}$. For all $k\geq 0$, by combining~\eqref{eq:NMLS} and~\eqref{eq:BLS} we have
\begin{equation}\label{eq:suff_dec_0}
\begin{aligned}
-\infty  < \varphi(x^{k+1}) = \varphi(\hat{x}^k + \lambda_k d^k) &   \leq  \varphi(\hat{x}^k)  - \alpha \lambda_k^2 \|d^k\|^2 \\
 & \leq  \varphi(x^k)  - a  \|d^k\|^2 -  \alpha \lambda_k^2 \|d^k\|^2    \\
 & \leq   \varphi(x^k)  - \eta  \left( 1 + \lambda_k^2 \right) \|d^k\|^2\\
 & =  \varphi(x^k)  - \eta \chi_k \|x^{k+1}-x^k\|^2,
\end{aligned}
\end{equation}
where $x^{k+1}-x^k= (1+\lambda_k) d^k$  and  $\chi_k :=  ( 1 + \lambda_k^2)/(1+\lambda_k)^2$  for all $k\geq 0$. 
In particular, observe that the first inequality in~\eqref{eq:suff_dec_0} holds as an equality when the optional linesearch in Step~6 of the algorithm is not performed, i.e., $\lambda_k=0$. By minimizing the value of $\chi_k$ with respect to $\lambda_k$ it is easy to see that
$\chi_k \geq  1/2$ for all $\lambda_k \geq 0$. Hence,~\eqref{eq:suff_dec_0} yields
\begin{equation}\label{eq:suff_dec}
-\infty  < \inf_{x\in\R^n} \varphi(x) \leq  \varphi(x^{k+1})  \leq   \varphi(x^k)  - \frac{\eta}{2}  \|x^{k+1}-x^k\|^2, \quad \forall k\geq 0,
\end{equation}
which in particular implies~\eqref{eq:dec1}. Then the sequence  $\bigl( \varphi(x^k)\bigr)_{k\in\N}$ is  nonincreasing and bounded from below. It follows that  $\bigl(\varphi(x^k)\bigr)_{k\in\N}$ converges almost surely to the random variable $\varphi^\ast_\infty : = \inf_{k\in\N} \varphi(x^k)$. Moreover, since $$|\varphi(x^k)| \leq \max \{|\varphi(x^0)|, |\inf_{x\in\R^n} \varphi(x)|\}, \text{ for all }k\geq0,$$ using Lebesgue's dominated convergence theorem, we get that  $\lim_{k\to\infty}  \Ebb\bigl[\varphi(x^k)\bigr] =  \Ebb\bigl[\varphi^\ast_\infty]$. This proves~\ref{it:t1-1}.

Now, by summing~\eqref{eq:dec1} for all $k\geq 0$, we get
\[
\sum_{k=0}^{\infty} \|x^{k+1}-x^k \|^2 \leq 2 \eta^{-1} \left( \varphi(x^0) -  \varphi^*_{\infty} \right).
\]
Then, by taking expectations in the above inequality we get that \eqref{eq:sum_deltak} holds. Finally, from~\eqref{eq:sum_deltak} we deduce that $\|x^{k+1}-x^k\|\toas 0$ and   $\Ebb[\|x^{k+1}-x^k\|]\to 0$, this concludes~\ref{it:t1-2} and the proof.
\end{proof}

Now, we turn  our attention to prove our main theorem, Theorem~\ref{t:2}, which determines that the  accumulation points of the sequence $(x^{k})_{k\in\N}$ are stationary points of problem~\eqref{eq:P1} almost surely. In order to obtain such a result, we need to assume that  the sequences  $(x^{k+1})_{k\in\N}$, $(d^{k})_{k\in\N}$ and $(\lambda_{k})_{k\in\N}$ are measurable with respect to the filtration generated by the random variables $i_0,\ldots, i_{k}$, that is, 
\begin{align*}
	\mathcal{F}_{k}:= \sigma\left(\left\{  i_0,\ldots, i_{k}\right\} \right).
\end{align*}

Observe that at every iteration the updates $x^{k+1}$ and $d^k$ are computed using only block-information of the objective function. Hence, they are not adequate to derive global stationary conditions for $\varphi$. This leads us to introduce the following auxiliary sequences.

First, for each $i\in{\{1,\ldots,N\}}$, let us denote as  $\bar{d}^{k,i}\in\R^n$ to the vector defined by Step~3 of Algorithm~\ref{alg:1} if $i_k = i$ and such that it provides a successful  direction for the  linesearch in \eqref{eq:NMLS}. Namely, $\bar{d}^{k,i}_{t}=0$ whenever $t\neq i$ and there exists $\tau_{k,i} >0$ (by~Proposition~\ref{prop:ls_wd})  such that 
\begin{equation}\label{eq:bard_ik}
\bar{d}^{k,i}_i \in \prox_{ \tau_{k,i} g_{i}}\left(   x^k_{i} - \tau_{k,i} \nabla_{i } f(x^k)  \right) - x^k_{i}
\end{equation}
and 
\begin{equation}\label{eq:decphi_ik}
\varphi(x^k+ \bar{d}^{k,i}) \leq  \varphi(x^k) - a\|\bar{d}^{k,i}\|^2.
\end{equation}
In particular, $\tau_{k,i}$ is the proximal stepsize accepted by the  linesearch of Algorithm~\ref{alg:1}  at iteration $k $ if  $i_k=i$. Further, we define the diagonal matrix $\varTh_k:=[\tau_{k,1} U_1,\ldots, \tau_{k,N} U_N]$, that contains the proximal stepsizes accepted in iteration $k$ for the different index selections. 

Similarly, let {$\varLam_k :=  [\lambda_{k,1} U_1,\ldots, \lambda_{k,N} U_N] $} be the diagonal matrix that  contains the stepsizes chosen  in  the $k$-th iteration   of Algorithm~\ref{alg:1}  by the additional linesearch  of Step~6   for the different index selections. Particularly, if this optional linesearch is not performed  at iteration $k$ with index $i_k=i$ then $\lambda_{k,i}=0$. Define the vectors
\begin{align}\label{defbard}
\bar{d}^k:= \sum_{i=1}^N \bar{d}^{k,i} \quad{ and } \quad \bar{\delta}^k := (I_n+\varLam_k) \bar{d}^k.
\end{align}
Thus, $ \bar{d}^k$ and $\bar{\delta}^k$ are random vectors $\mathcal{F}_{k}$-measurable.
Observe that, by~\eqref{eq:bard_ik} and variable separability in $g$, the vector $\bar{d}^k$  admits the  compact representation
\begin{equation}\label{eq:bard_prox}
\bar{d}^k\in \argmin_{d\in\R^n} \left\{ g(x^k+d) + \langle \nabla f(x^k),d\rangle + \frac{1}{2}\langle d, \varTh^{-1}_k d\rangle\right\}.
\end{equation}
The sequence $(\bar{d}^k)_{k \in \N}$, defined in~\eqref{defbard}, will play a key role in linking the sequence $(x^k)_{k \in \N}$ to a stationary point of problem~\eqref{eq:P1}. Before establishing this connection, we introduce two technical propositions. 

{
The first of them establishes that, under Assumption~\ref{assumption:1}, the sequence of matrices $\varTh_k$ associated to a  bounded sequence generated by Algorithm~\ref{alg:1} are uniformly positive definite. More precisely, we state the following proposition.

\begin{proposition}\label{prop:tau_bound}
	Suppose Assumption~\ref{assumption:1} holds. Let  $(x^{k_j})_{j\in\N}$ be a bounded subsequence of  the sequence $(x^k)_{k\in\N}$ generated by Algorithm~\ref{alg:1}. Then there exists a positive constant $\munderbar{\theta}>0$ such that $\varTh_{k_j} \succeq \munderbar{\theta} I_n$ for all $j\geq 0$.
\end{proposition}

\begin{proof}
	By way of contradiction, we assume that the result does not hold. By the diagonal structure of $\varTh_{k_j}$, this implies that for at least one index $\iota\in\{1,\ldots,N\}$ there exists a subsequence $(\tau_{k_\nu,\iota})_{\nu\in\N}$ of $  (\tau_{k_j,\iota})_{j\in\N}$ such that  $\tau_{k_\nu,\iota} \to 0^+$.

	Due to the latter, we can further assume that $\tau_{k_{\nu},\iota} < \munderbar{\tau}$, for all $\nu\in\N$. Recall that, according to the backtracking procedure, $\tau_{k_v,\iota}$ takes the form $\tau_{k_\nu,\iota} = \beta^{t_\nu} \bar{\tau}_{k_\nu}$ for some $t_{\nu}\in\N{\setminus}\{0\}$.\footnote{Rigorously, the counter $t_\nu$ should also be written with subscript $\iota$ (i.e., $t_{\nu,\iota}$) as $t_\nu$ may vary depending on the coordinate block that is being updated. However, we omit it to simplify notation. Similarly occurs with $\bar{\tau}_{k_\nu}$.}   Hence, the stepsize $\beta^{t_{\nu}-1}\bar{\tau}_{k_\nu}$  was not accepted.
    
    Now, let $\hat{d}^{k_\nu,\iota}$ be given by the proximal step~\eqref{eq:bard_ik} with  stepsize $\beta^{t_{\nu}-1}\bar{\tau}_{k_\nu}$, i.e., $\hat{d}^{k_\nu,\iota}_{i}=0$ if $i\neq \iota$ and
	\begin{equation}\label{eq:prox_dhat}
		\begin{aligned}
			\hat{d}^{k_{\nu},\iota}_{\iota} \in \argmin_{s} \left\{ g_{\iota}(x^{k_\nu}_{\iota}+s) + \langle \nabla_{\iota} f (x^{k_\nu}), s\rangle + \frac{1}{2\beta^{t_{\nu}-1} \bar{\tau}_{k_\nu}}\|s\|^2 \right\},  \quad \forall \nu\geq0.
		\end{aligned}
	\end{equation}
	On the one hand, since the stepsize $\beta^{t_{\nu}-1}\bar\tau_{k_\nu}$ was not accepted, then $\|\hat{d}^{k_\nu,\iota}\|\neq 0$ (by Step~4 of Algorithm~\ref{alg:1})   and  condition~\eqref{eq:decphi_ik} failed. The latter implies that
	\begin{equation}\label{eq:F_hat}
		\varphi(x^{k_\nu} + \hat{d}^{k_\nu,\iota})  > \varphi(x^{k_\nu}) - a \|\hat{d}^{k_\nu,\iota}\|^2, \quad \forall \nu\geq 0.
	\end{equation}
	On the other hand, \eqref{eq:prox_dhat} yields
	\[
	g(x^{k_\nu} + \hat{d}^{k_\nu,\iota}) +\langle \nabla_{\iota} f (x^{k_\nu}), \hat{d}^{k_{\nu},\iota}_{\iota} \rangle + \frac{1}{2\beta^{t_{\nu}-1}\bar{\tau}_{k_\nu}}\|\hat{d}^{k_{\nu},\iota}_{\iota}\|^2 \leq  g(x^{k_\nu} ).
	\] 
	Subtracting the above  inequality to~\eqref{eq:F_hat} we obtain
	\begin{equation}\label{eq:tau_f1}
		f(x^{k_\nu} + \hat{d}^{k_\nu,\iota})  > f(x^{k_\nu}) + \langle \nabla_{\iota f}(x^{k_\nu}), \hat{d}^{k_{\nu},\iota}_{\iota} \rangle +  \left( \frac{1}{2\beta^{t_\nu-1}\bar{\tau}_{k_\nu}} - a\right) \|\hat{d}^{k_\nu,\iota}\|^2, \quad \forall\nu\geq 0.
	\end{equation}
	
	Now, using the fact that  $(x^{k_\nu})_{\nu\in\N}$ is bounded, by possible resorting to a smaller subsequence, we  can assume without loss of generality  that $(x^{k_{\nu}})_{\nu\in\N}$ converges to some accumulation  point  $x^*$. Particularly, due to the fact that $\varphi(x^k)$ is monotonically nonincreasing (see Theorem \ref{t:1}) and lower semicontinuity, we have that  $x^*\in\dom g$. Hence, by noting  that $\beta^{t_{\nu}-1}\bar{\tau}_{k_\nu} = \beta^{-1}\tau_{k_\nu,\iota}$,  we have that 
    \[
    x^{k_\nu}_{\iota} + \hat{d}^{k_\nu,\iota}_{\iota} \in \prox_{\beta^{-1} \tau_{k_\nu,\iota} g_{\iota}} \bigl( x^{k_\nu}_{\iota}  - \beta^{-1}\tau_{k_\nu,\iota} \nabla_{\iota}  f (x^{k_\nu} )  \bigr).
    \]
    {Since $ x^{k_\nu}_{\iota}  - \beta^{-1}\tau_{k_\nu,\iota} \nabla_{\iota}  f (x^{k_\nu} ) \to x^\ast_{\iota}$,  Lemma~\ref{lemma:proxcoledk}  yields that $x^{k_\nu}_\iota + \hat{d}^{k_\nu,\iota}_\iota \to x^*_{\iota}$ as $\nu\to\infty$. Consequently, $x^{k_\nu} + \hat{d}^{k_\nu,\iota} \to x^*$ as $\nu\to\infty$, since by construction $\hat{d}^{k_\nu,\iota}_i =0$ for $i\neq \iota$.} Furthermore, by Assumption~\ref{assumption:1} there exists a positive constant~$L_{\iota}(x^*)$ such that   for all sufficiently large $\nu\in\N$ it holds
	\begin{equation}\label{eq:tau_f2}
		f(x^{k_\nu} + \hat{d}^{k_\nu,\iota})  \leq  f(x^{k_\nu}) + \langle \nabla_{\iota} f(x^{k_\nu}), \hat{d}^{k_{\nu},\iota}_{\iota} \rangle +  \frac{L_{\iota}(x^*)}{2} \|\hat{d}^{k_\nu,\iota}\|^2.
	\end{equation}
	Putting together~\eqref{eq:tau_f1} and~\eqref{eq:tau_f2} and dividing by $\|\hat{d}^{k_\nu,\iota}\|^2$  we get
	\[
	\frac{\beta}{2\tau_{k_\nu,\iota}} - a<   \frac{L_{\iota}(x^*)}{2}.
	\]
	As the left hand-side  of the above inequality tends to $+\infty$ as $\nu\to\infty$, this leads to the desired contradiction.
\end{proof}
}

{
\begin{remark}
Some comments regarding Proposition~\ref{prop:tau_bound} follow.
\begin{enumerate}[label=(\roman*)]
\item Under the standard, more restrictive  assumption,   of $f$  having a block-wise (globally) Lipschitz continuous gradient, the boundedness assumption in Proposition~\ref{prop:tau_bound}   becomes unnecessary. In this case, $L_{i_k}:= L_{i_k}(x^k)$ for all $k\in\N$. Consequently, \eqref{eq:sat_ls2} directly yields a  uniform lower  bound for the stepsize, given by
\[
\min\left\{\munderbar{\tau}, \frac{\beta}{  \max_{i\in{\{1,\ldots,N\}}} L_i + 2a  }\right\} \leq \tau_k, \quad \forall k\in\N,
\]
where we recall that $\beta$ is the backtracking parameter. In contrast, if the gradient of $f$  is merely block-wise locally Lipschitz  continuous,  the local Lipschitz constants can grow arbitrarily large along an unbounded sequence $\{x^k\}_{k\in\N}$, leading to the stepsizes vanishing to zero. 

\item As just mentioned,  Proposition~\ref{prop:tau_bound} (and also our main convergence result; see Theorem~\ref{t:2} below) assumes the existence of a bounded subsequence of the iterates generated by Algorithm~\ref{alg:1}. Recall that this is a standard assumption required for ensuring (subsequential) iterate convergence of optimization algorithms for nonconvex problems, even under global Lipschitz conditions (see, e.g.,\cite[Theorem~1]{MR3293507}) or the well-known Kurdyka--{\L}ojasiewicz framework (see, e.g.,~\cite[Theorem~4.1]{MR3572363}). This boundedness assumption immediately  holds for problems where the lower level sets of the objective function $\varphi$ are bounded or, similarly, whenever the domain of the prox-bounded function $g$ is bounded.
\end{enumerate}
\end{remark}
}

The following result relates the expected norms of $\bar{\delta}^k$ and $\bar{d}^k$ to the expected norm of $x^{k+1} - x^k$ for all $k \geq 0$. In particular, it shows that the sequence $(\Ebb[\|\bar{d}^k\|^2])_{k \in \N}$ is summable. This summability will be essential for our main theorem, as it allows us to take the limit in the proximal step~\eqref{eq:bard_prox} along a suitable subsequence.


\begin{proposition}\label{prop:exp_bar}
	Suppose Assumption~\ref{assumption:1} holds. Given $x^0\in\dom \varphi$, let $(x^k)_{k\in\N}$ be the sequence generated by Algorithm~\ref{alg:1}. In addition, assume that $\bigl(x^{k+1}, d^{k},\lambda_k\bigr)$ is $\mathcal{F}_{k}$-measurable for all $k\geq 0$.  Then the sequences $(\bar{\delta}^k)_{k\in\N}$ and $(\bar{d}^k)_{k\in\N}$, defined in \eqref{defbard}, satisfy that
	\begin{equation}\label{eq:t2-claim1}
		\Ebb\bigl[  \|x^{k+1}-x^k \|^2 \bigr]  \geq  p_{\min}\Ebb \bigl[ \| \bar{\delta}^{k} \|^2 \bigr] \geq  p_{\min}  \Ebb \bigl[  \| \bar{d}^k\|^2 \bigr] , \quad \forall k\geq 0,
	\end{equation}
	where $p_{\min} := \min\{p_i \, : \, i=1,\ldots N\}$. Consequently,
	\begin{equation}\label{eq:t2-claim2}
		\sum_{k=0}^{\infty   }      \Ebb \bigl[ \| \bar{\delta}^{k} \|^2 \bigr]   <  \infty,  \quad 	\sum_{k=0}^{\infty   }      \Ebb \bigl[ \| \bar{d}^{k} \|^2 \bigr]   < \infty, \quad 
 \bar{\delta}^{k}  \toas  0 \quad \text{and} \quad  \bar{d}^{k} \toas  0.
	\end{equation} 
     	 Furthermore, there exists a constant $C\geq 0$ such that 
	 \begin{equation}\label{eq:it:t3-1}
	 	\min_{0\leq k \leq \ell } \mathbb{E}\bigl[ \| \bar{\delta}^{k} \|^2 \bigr]  \leq   \frac{C}{ \ell +1 }, \quad \forall \ell \geq 0.  
	 \end{equation}
	\end{proposition}
	\begin{proof}
		Indeed, for any $k\geq 0$, {using the measurability of $(x^{k+1},d^k,\lambda_k)$ and since $\mathbb{P}(i=i_k | \mathcal{F}_{k-1}) = p_i$}, we have
		\begin{equation}\label{eq:p_min}
		\begin{aligned}
			\Ebb\bigl[ \|x^{k+1}-x^k\|^2 \, | \, \mathcal{F}_{k-1} \bigr]      
            & = \sum_{i=1}^N p_i (1+\lambda_{k,i})^2 \|\bar{d}^{k,i}\|^2\\ 
			&  \geq  p_{\min} \sum_{i=1}^N   (1+\lambda_{k,i})^2\|\bar{d}^{k,i}\|^2= p_{\min} \| \bar{\delta}^k\|^2 \geq p_{\min} \| \bar{d}^k\|^2 ,
		\end{aligned}
		\end{equation}
		with the convention $\mathcal{F}_{-1} := \{ \emptyset, \Omega\}$.   Hence,
		\begin{equation*}
         \begin{aligned}
			\Ebb\bigl[  \|x^{k+1}-x^k \|^2 \bigr] &  = \Ebb \biggr[   \Ebb \bigl[ \|x^{k+1}-x^k\|^2 \, | \, \mathcal{F}_{k-1} \bigr] \biggl] \\
             & \geq \,    p_{\min}\Ebb \bigl[ \| \bar{\delta}^{k} \|^2 \bigr] \geq  p_{\min}\Ebb \bigl[ \| \bar{d}^{k} \|^2 \bigr] , \quad \forall k\geq 0.
            \end{aligned}
		\end{equation*}
			Furthermore, the above equation and Theorem \ref{t:1}\ref{it:t1-2} imply that the series in  \eqref{eq:t2-claim2} converge.
			
			Now, let us note that by exchanging summation and expectation (see, e.g., \cite[Corollary 8.2.11]{MR4704094}), we obtain  that the convergence of the series in \eqref{eq:t2-claim2} is equivalent to 
			\begin{equation*} 
				\Ebb \left[	\sum_{k=0}^{\infty}      \| \bar{\delta}^{k} \|^2 \right]   <  \infty  \quad \text{and}  \quad	  \Ebb \left[ \sum_{k=0}^{\infty}    \| \bar{d}^{k} \|^2 \right]   < \infty.
			\end{equation*}
			This implies that, almost surely 
			\begin{equation*} 
				\sum_{k=0}^{\infty}      \| \bar{\delta}^{k} \|^2   <  \infty  \quad \text{and}  \quad	   \sum_{k=0}^{\infty}    \| \bar{d}^{k} \|^2    < \infty,
			\end{equation*}
		  which in turn ensures that the sequences $(\bar{d}^{k})_{k\in\N}$ and $( \bar{\delta}^{k})_{k \in \N}$ almost surely converge to zero, and this completes the proof of \eqref{eq:t2-claim2}.

		  Finally, by taking expectations in~\eqref{eq:dec1} we obtain
		  \[
		  \Ebb[ \varphi(x^{k+1} )] \leq \Ebb[\varphi(x^k) ]- \frac{\eta}{2} \Ebb [\|x^{k+1}-x^k\|^2 ], \quad \forall k\geq 0.
		  \]
		  Summing the above inequality over $k =0,\ldots, \ell$ and telescoping yields
		  \begin{equation}\label{eq:t3-1}
		  	\frac{\eta}{2}\sum_{k=0}^\ell \Ebb [\|x^{k+1}-x^k\|^2 ] \leq \Ebb[ \varphi(x^0)] - \Ebb [ \varphi(x^{\ell+1}) ]  \leq    \varphi(x^0) - \Ebb [ \varphi^*_{\infty} ],
		  \end{equation}
		  where $\varphi^*_{\infty} := \inf_{k\in \N} \varphi(x^k)$ is defined as in Theorem~\ref{t:1}. Now,  \eqref{eq:t2-claim1} and \eqref{eq:t3-1} give
		  \[
		  p_{\min}\min_{0\leq k\leq \ell}  \Ebb [\|\bar{\delta}^k\|^2 ]  \leq  \min_{0\leq k\leq \ell}  \Ebb [\|x^{k+1}-x^k\|^2 ] \leq    \frac{2}{ \eta (\ell+1) } \left( \varphi(x^0) - \Ebb [ \varphi^*_{\infty} ] \right),
		  \]
		  which proves~\eqref{eq:it:t3-1}.

		\end{proof}
		
		\begin{remark}\label{remarkasc}

First, let us observe that the $\mathcal{F}_{k-1}$-measurability of $x^k$ is essential in our analysis. Note that if any of the sequences generated by the algorithm fails to form an adapted process with respect to the filtration $(\mathcal{F}_k)_{k\in\mathbb{N}}$, then the first inequality in \eqref{eq:p_min} may no longer hold.

		Second, it is well known that convergence in expectation does not necessarily imply almost sure convergence (see, e.g., \cite[Example 1.5.5]{MR4704094}). Therefore, from \eqref{eq:t2-claim1} and the fact that $\mathbb{E}\bigl[\|x^{k+1} - x^k\|^2\bigr] \to 0$, we cannot conclude that the random vector $\bar{d}_k$ converges to zero almost surely. Consequently, we cannot directly pass to the limit in \eqref{eq:bard_prox} to guarantee that {accumulation points} of the sequence $(x^k)_{k\in\mathbb{N}}$ are stationary points of problem \eqref{eq:P1}. To overcome this limitation, the summability of the sequences $\bigl(\mathbb{E}[\|\bar{d}^k\|^2]\bigr)_{k\in\mathbb{N}}$ and $\bigl(\mathbb{E}[\|\bar{\delta}^k\|^2]\bigr)_{k\in\mathbb{N}}$, as established in \eqref{eq:t2-claim2}, plays a crucial role. This summability enables us to establish the almost sure convergence of both $(\bar{d}^k)_{k\in\mathbb{N}}$ and $(\bar{\delta}^k)_{k\in\mathbb{N}}$ to zero, and then, in our main theorem the argument described above can be applied. 
			
		\end{remark}

We are now in position to present our main result.

\begin{theorem}\label{t:2}
Suppose  Assumption~\ref{assumption:1} holds. Given $x^0\in\dom \varphi$, let $(x^k)_{k\in\N}$ be the sequence generated by Algorithm~\ref{alg:1}. Assume that  $\bigl(x^{k+1}, d^{k},\lambda_k\bigr)$ is $\mathcal{F}_{k}$-measurable for all $k\geq 0$. Then, almost surely, the following statements hold:
\begin{enumerate}[label=(\roman*)] 
\item\label{it:t2-1} Any accumulation point $x^*$ of the sequence $(x^{k})_{k\in\N}$ is a stationary point of  problem~\eqref{eq:P1}. 

\item\label{it:t2-2}  If $(x^k)_{k\in\N}$ has at least one isolated accumulation point, then the whole sequence  $(x^k)_{k\in\N}$ converges to a stationary point of~\eqref{eq:P1}.

\end{enumerate}
\end{theorem}

\begin{proof} 
\ref{it:t2-1} Let $(x^{k_\nu})_{\nu\in\N}$ be a subsequence of $(x^k)_{k\in\N}$ converging to the accumulation point $x^*$.
Since the sequence  $(x^{k_\nu})_{\nu\in\N}$ is bounded, by Proposition~\ref{prop:tau_bound} there exists $\munderbar{\theta}>0$ such that $\munderbar{\theta} I_n \preceq \varTh_{k_\nu}   \preceq \min{\{\bar{\tau},\tau^{g}\}} I_n$ for all $\nu\in\N$, where $\tau^{g}:= \min\{\tau^{g_i} \, : \, i=1,\ldots, N\}$. Thus we can further assume that  $\varTh_{k_\nu}\to\varTh^*\succeq\munderbar{\theta}I_n$, with $\varTh^*:=  [\tau_{1}^* U_1,\ldots, \tau_{N}^* U_N]$ and $\tau_i^\ast >0$ for all $i\in{\{1, \ldots,N\}}$. 
 
Observe that replacing $k$ by $k_\nu$ in~\eqref{eq:bard_prox} yields 
\begin{equation}\label{eq:graph_c}
\bigl( x^{k_\nu} - \varTh_{k_\nu} \nabla f(x^{k_\nu}),  x^{k_\nu} + \bar{d}^{k_\nu} \bigr)_{\nu\in\N} \subset \mathrm{graph} \prox_{\varTh_{k_\nu} g},
\end{equation}
where $\prox_{\varTh_{k_\nu} g}$ is defined as the proximity operator of the function $\sum_{i=1}^N \tau_{k_\nu,i} \, g_i$.
Due to variable separability,~\eqref{eq:graph_c} is in turn equivalent to
\begin{equation}\label{eq:graph_nc}
\bigl(x^{k_\nu}_i - \tau_{k_\nu,i} \nabla_i f(x^{k_\nu}), x^{k_\nu}_i + \bar{d}^{k_\nu,i}\bigr)_{\nu\in\N} \subset \mathrm{graph} \prox_{\tau_{k_\nu,i} g_i}, \quad  i=1,\ldots,N.
\end{equation}

Now, note that by Proposition~\ref{prop:exp_bar}   it holds that $\bar{d}^k\toas 0$.
Hence, taking the limit as $\nu\to\infty$ in the sequences in the left-hand side of \eqref{eq:graph_nc} yields
\begin{equation}\label{eq:toascp}
\bigl(x^{k_\nu}_i - \tau_{k_\nu,i} \nabla_i f(x^{k_\nu}), x^{k_\nu}_i + \bar{d}^{k_\nu,i}\bigr) \toas (x^*_i - \tau_i^* \nabla_i f(x^*), x^*_i ),  \quad  i=1,\ldots,N.
\end{equation}

Furthermore,  {by~\eqref{eq:graph_nc}-\eqref{eq:toascp} and $\tau_{k_\nu,i}\to\tau^*_i$,  the last assertion in~\cite[Theorem~1.25]{MR1491362} ensures} that,  for all $i\in\{1,\ldots,N\}$,
\[
 \quad x^*_i \in \prox_{\tau_i^* g_i} \bigl(x_i^* -\tau_i^* \nabla_i f(x^*) \bigr)  \text{ a.s.},
\]
which again due to variable separability can be expressed in compact form   as
\begin{equation}\label{eq:prox_f}
x^* \in \prox_{\varTh^* g} \bigl(x^* - \varTh^* \nabla f(x^*) \bigr) \text{ a.s.}
\end{equation}
Since the first order necessary optimality condition of~\eqref{eq:prox_f} leads to
\[
0 \in \partial g(x^*) + (\varTh^*)^{-1} \bigl(x^* - \bigl(x^* - \varTh^* \nabla f(x^*)\bigr) \bigr) \Longrightarrow  0 \in \partial g(x^*) + \nabla f(x^*),
\] 
this concludes that $x^*$ is a stationary point of~\eqref{eq:P1} almost surely.


\ref{it:t2-2}  Let $x^*$ be an isolated accumulation point of the sequence $(x^k)_{k\in\N}$. Since, by Theorem~\ref{t:1}\ref{it:t1-2}, the sequence $(x^k)_{k\in\N}$ verifies the so-called \emph{Ostrowski property} almost surely, namely,
\[
\lim_{k\to\infty} \|x^{k+1}-x^k\| = 0 \text{ a.s.,}
\]
then $x^k\toas x^*$ as $k\to\infty$, by~\cite[Proposition~8.3.10]{MR1955649}. Finally, by~\ref{it:t2-1} we deduce that $x^*$ is a stationary point of~\eqref{eq:P1} almost surely, and that concludes the proof of the result.
\end{proof}



\section{Numerical experiments}\label{sect:numerical}

In this section, we test numerically the performance of different versions of the Adaptive Randomized Block Proximal Gradient (ARBPG) method in Algorithm~\ref{alg:1}. The experiments were run on a computer of AMD Ryzen 9 7950X 16-Core Processor, 4.50 GHz with 64GB RAM, under Windows 10 (64-bit). {The organization is as follows. In Section~\ref{sect:num1},  we develop a concrete implementation of the boosted linesearch of our algorithm.   In Section~\ref{sect:num2}, we test the performance of the basic version (without boosted linesearch) of Algorithm~\ref{alg:1} in  experiments in nonnegative matrix factorization. More precisely, we analyze the effect of different choices of block sizes and compare against the recognized \emph{inertial proximal alternating linearized minimization} method~\cite{MR3572363}. Finally, in Section~\ref{sect:num3}, we illustrate the benefit of the boosted linesearch for problems where the Lipschitz moduli cannot be approximated accurately.} The Python source code and data for the
experiments are available and can be consulted  at \href{https://github.com/DavidTBelen/ARBPG}{https://github.com/DavidTBelen/ARBPG}.

\subsection{Algorithm implementation}\label{sect:num1}

We distinguish two variants of Algorithm~\ref{alg:1} depending on whether  the optional linesearch in Step~6 of the method is activated. When such a linesearch is not performed, for instance by setting $\lambda_k=0$ for all $k\geq0$,  we denote the resulting scheme as ARBPG. On the other hand,  we use   ARBPG-B for  the version of Algorithm~\ref{alg:1} that incorporates a  \emph{boosted linesearch}~\cite{BDSA2025,MR4078808,MR3785672} in  Step~6.

As commented before, in the work~\cite{BDSA2025} the authors introduced a procedure for efficiently computing boosted linesearches with a fixed number of steps when the direction defined by the proximity operator is not a descent direction.
Here, we present a new approach  that does not require to predetermine a fixed number of attempts in the linesearch.  Specifically, the method ARBPG-B considered here makes use of Algorithm~\ref{alg:bls} to adjust the parameter $\lambda_k$ in Step~6 of Algorithm~\ref{alg:1}. The main characteristic of the procedure detailed in Algorithm~\ref{alg:bls} is that, instead of performing the backtracking linesearch in the interval $\{\hat{x}^k + \lambda d^k :  \lambda > 0\}$, it considers the interval $\{ x^k + \varsigma d^k  :  \varsigma >0\}$ and initializes the search at some  step $\varsigma_k>1$. This way, if after some backtracking steps condition~\eqref{eq:BLS} did not hold and  $\varsigma_k $ takes a value smaller than 1, it simply stops the linesearch by setting  $\varsigma_k=1$, in which case the resulting point is $\hat{x}^k$ (same as when $\lambda$ takes value 0), and~\eqref{eq:BLS} holds trivially as an equality. 
\begin{algorithm}[h!]
	\caption{Boosted linesearch procedure for Step~6 in Algorithm~\ref{alg:1} without differentiability of $g$}\label{alg:bls}
	\begin{algorithmic}[1]
		\Require{ $x^k$; $\hat{x}^k$; $d^k$; $\alpha>0$ and $\rho\in{]0,1[}$.}
		\State{Choose any $\overline{\varsigma}_k\geq1$. Take $l=0$ and set $\varsigma_k:=\rho^l\overline{\varsigma}_k$.}
		\State{\textbf{while} $\varsigma_k>1$ and 
		\begin{align*}
		\varphi(x^k+\varsigma_k d^k)  > \varphi(\hat{x}^k) - \alpha (\varsigma_k-1)^2 \|d^k\|^2
		\end{align*}    
		\hspace{\algorithmicindent}\textbf{do} $l\leftarrow l+1$, $\varsigma_k:=\rho^l \varsigma_k$.} 
		\State{\textbf{if} $\varsigma_k>1$  \textbf{then} }\\
		\hspace{\algorithmicindent} \Return{$\lambda_k := \varsigma_k-1$;} 
		\State{\textbf{else}}\\
		\hspace{\algorithmicindent} \Return{$\lambda_k := 0$.} 
		\State{\textbf{end if}}
	\end{algorithmic}
\end{algorithm}

{In all of our experiments and for both ARBPG and ARBPG-B, we set $\munderbar{\tau} = 10^{-8}$, $\bar{\tau} = 10^8$. Recall that these parameters are only required for theoretical validation and do not have any practical impact on the algorithm performance. When needed, the parameter for the Armijo-type linesearch is taken as $a=10^{-4}$, which is a value commonly used in practice, and for the backtracking parameter we set  $\beta=0.9$. According to Remark~\ref{remark:trial_step}, the choice of trial stepsize $\bar{\tau}_k$ in Step~2 of ARBPG will try to approximate the reciprocal of the local Lipschitz moduli. Specific choices shall be given below for each problem.}

For ARBPG-B we set 
$\alpha=0.1$. Furthermore, the boosted linesearch at each iteration of ARBPG-B also requires the selection of a backtracking parameter $\rho$ and an initial trial stepsize $\overline{\varsigma}_k$ for Algorithm~\ref{alg:bls}. Particularly, we set $\rho=0.5$ while~$\overline{\varsigma}_k$ is set adaptively at every iteration following  a slightly modified version of the \emph{self-adaptive trial stepsize} presented in~\cite[Algorithm~2]{BDSA2025}, which is in turn based on the one proposed in~\cite{MR4078808}. Specifically, we employ the adaptive procedure gathered in Algorithm~\ref{alg:fortracking} with the choices $\bar{\varsigma}_0=3$ and $\delta=2$. 
\begin{algorithm}[H]
	\caption{Self-adaptive trial stepsize}\label{alg:fortracking}
	\begin{algorithmic}[1]
		\Require{$\delta>1$ and $\overline{\varsigma}_0> 1$. Obtain $\varsigma_k$ from $\overline{\varsigma}_k$ by Steps $1$-$2$ of  Algorithm~\ref{alg:bls}.}
		\If{$\varsigma_k = \bar{\varsigma}_k$}
		\State{set $\overline{\varsigma}_{k+1}:=\delta\overline{\varsigma}_k$;}
		\Else
		\State{set $\overline{\varsigma}_{k+1}:=\max{\{ \overline{\varsigma}_0, \varsigma_k}\}$.}
		\EndIf
	\end{algorithmic}
\end{algorithm}
 The interpretation of the adaptive procedure in Algorithm~\ref{alg:fortracking} is simple. For any $k\in\N$, the starting stepsize $\overline{\varsigma}_{k+1}$ for the next iteration of the boosted linesearch is determined depending on the performance of the current boosted linesearch. If at the current iteration the boosted linesearch succeeded at the first attempt (i.e., when $\varsigma_k =\overline{\varsigma}_k$), then the trial stepsize for the next iteration is increased by setting $\overline{\varsigma}_{k+1} := \delta \varsigma_k$, with $\delta >1$. Otherwise, $\overline{\varsigma}_{k+1}$ is set as the maximum between the default initial stepsize and the smallest stepsize tested in the previous iteration, namely, $\overline{\varsigma}_{k+1}:=\max\{\overline{\varsigma}_0, \varsigma_k \}$. Observe that  at iteration $k$ of Algorithm~\ref{alg:1} the choice of $\lambda_k$ resulting from Algorithms~\ref{alg:bls} and~\ref{alg:fortracking} is $\mathcal{F}_{k}$-measurable, so the assumptions of Proposition~\ref{prop:tau_bound} and Theorem~\ref{t:2} hold.

{Table~\ref{t:algorithm_parameters} summarizes the parameter configuration for ARBPG and ARBPG-B used for all the experiments in this section. Finally, we emphasize that although the optimal choice of the boosted parameters is problem-dependent, the parameter setting used here mimics the one in~\cite{BDSA2025}  which has proven to consistently outperform non-boosted variants of algorithms across different problems. As an alternative, one could conduct a thorough study to determine optimal parameter configurations. However, we prefer to use this standard configuration to show that, even without an optimal parameter selection for the specific problem, the boosted linesearch may lead to advantages in algorithmic performance.}

\begin{table}[h!]
\centering
\caption{Summary of the parameters used by  ARBPG and ARBPG-B.}
\label{t:algorithm_parameters}
\begin{tblr}{
  width = \linewidth,
  colsep = 4pt,          
  rows = {abovesep=5pt, belowsep=5pt},
  colspec = {
    |Q[c, valign=m] |      
    X[c, valign=m] |      
    X[l, valign=m] |      
  },
}
\hline
Instances of Algorithm~\ref{alg:1} & {Adaptive linesearch} & {Boosted linesearch} \\
\hline
ARBPG   & \SetCell[r=2]{c} $\underline{\tau}=10^{-8}$, $\bar{\tau}=10^8$,\newline $a=10^{-4}$, $\beta= 0.9$ 
                         & $\lambda_k=0$ for all $k \in \mathbb{N}$ \\
\hline
ARBPG-B &                & Algorithms~\ref{alg:bls} and~\ref{alg:fortracking} with\newline $\alpha=0.1$, $\rho=0.5$, $\overline{\varsigma}_0=3$, $\delta = 2$ \\
\hline
\end{tblr}
\end{table}

\subsection{Nonnegative matrix factorization}\label{sect:num2}
The \emph{Nonnegative Matrix Factorization} (NMF) problem consists in approximating a matrix $A\in\R^{m\times n}$ as a product of two low rank matrices $\mU\in\R^{p\times m}$ and $\mV\in\R^{p\times n}$ with nonnegative entries. NMF finds applications in multiple fields, such as text mining~\cite{MR2388467}, face recognition~\cite{LeeNMF} or spectral data analysis~\cite{PAUCA200629}, among other areas. Mathematically, the problem can be formulated as 
\begin{equation}\label{eq:NMF1}
\min_{\mU, \mV} \quad \frac{1}{2}\| A- \mU^T \mV  \|_F^2  \quad \text{s.t.} \quad \mU\in\R_+^{p\times m}, \; \mV\in\R_+^{p\times n},
\end{equation}
where $\|\cdot\|_F$ denotes  the \emph{Frobenius norm}.  This problem can  be posed in  the form of \eqref{eq:P1} by setting
\begin{equation}\label{eq:NMF2}
f(\mU,\mV) := \frac{1}{2}\| A- \mU^T\mV\|_F^2 \quad \text{and}  \quad  g(\mU,\mV):=  \iota_{\R_+^{p\times m}}(\mU) + \iota_{\R^{p\times n}_+}(\mV),
\end{equation}
with  $\iota_C$ denoting  the \emph{indicator function} of a set $C$.  Observe that the partial derivatives of $f$ with respect to $\mU$ and $\mV$ are given by 
\begin{equation}\label{eq:NMF_grad}
\nabla_{\mU} f(\mU,\mV) = - \mV(A^T -\mV^T\mU) \quad \text{and}\quad  \nabla_{\mV} f(\mU,\mV) = -\mU (A - \mU^T \mV),
\end{equation}
respectively. It is easy to see that  they are locally Lipschitz continuous, but not globally Lipschitz. In addition, $g$ can clearly be defined entrywise and its proximity operator is computed by projecting element-wise onto the nonnegative real semi-line $\R_+$. Thus the problem perfectly fits into the framework delineated by Assumption~\ref{assumption:1}.

Block-coordinate methods are among the preferable choices for dealing with nonnegative matrix factorization problems in the form of~\eqref{eq:NMF1}. We refer the reader to~\cite{MR3158055} for a review of some strategies for splitting the variables of the problem. {Perhaps, the most popular method used for NMF is the one known as \emph{Proximal Alternating Linearized Minimization} (PALM), proposed by Bolte, Sabach and Teboulle~\cite{MR3232623}; or its inertial version iPALM, devised by Pock and Sabach~\cite{MR3572363}. A stochastic extension of PALM has also been  proposed in~\cite{MR4354998}. These methods tackle~\eqref{eq:NMF1} by alternating between proximal problems  involving the variables $\mU$ and $\mV$ independently. A sufficient decrease condition is guaranteed when the stepsizes for such proximal steps are chosen to be smaller than the reciprocal of the  Lipschitz constants of the partial gradients.  This approach makes use of the fact that these partial gradients  are global once the  variable held constant in the derivation is fixed. Indeed, as explained in~\cite{MR3572363},  the Lipschitz moduli of $\nabla_V f$ and $\nabla_W f$ are directly deduced from the linearity of~\eqref{eq:NMF_grad}, and result in
\[
L_V(W) = \|\mV\mV^T||_2 \quad \text{and} \quad L_W(V) = \|\mU\mU^T||_2,
\]
where $\|\cdot\|_2$ stands for the matrix $2$-norm, i.e., the largest singular value.  

 In order to implement ARBPG we consider blocks defined by the rows of $\mU$ and $\mV$. The Lipschitz moduli of the partial gradients given by the blocks are then obtained as
\[
L_{\mU_i}(\mV) = \|\mV_i \mV_i^T||_2 \quad \text{and} \quad L_{\mV_i}(\mU) = \|\mU_i\mU_i^T||_2,
\]
where $\mU_i$ and $\mV_i$ are the submatrices of $\mU$ and $\mV$ containing only the rows associated with the $i$-th block. Similarly to what happens with PALM, if the trial stepsize $\bar{\tau}_k$ in Step~3 of ARBPG is taken smaller than the reciprocal of the Lipschitz modulus, then the Armijo-type linesearch is accepted in the first iteration. To establish a fair comparison between the two methods, in the following experiment we will always select the corresponding stepsize by multiplying the reciprocal of the Lipschitz constant by 0.95. In addition, as in this setting we assume that the exact Lipschitz moduli can be computed efficiently, the boosted linesearch is not expected to provide any significant impact. Hence, ARBPG-B is ruled out of the tests in this section.

Note that in this context, if the blocks of ARBPG are taken to be the two matrices $\mU$ and $\mV$, then  a randomized version of PALM is recovered. Hence, it makes sense to study whether, by taking smaller blocks, ARBPG obtains an improved performance. This is what we investigate in our first experiment in the next section.}

\subsubsection{Experiment on image compression}

 We aim to  compress $4$ color images   of Chilean  landscapes with different resolutions.  This can be achieved through NMF by choosing a rank $p$  satisfying the  condition $(m+n) \times  p < m \times n$, where $m \times n$ is the pixel size of the original image.  Hence the total number of variables in the problem is $p \times (m+n)$. Table~\ref{tab:images} summarizes the resolution of each  image considered and the rank taken for the nonnegative matrix factorization. 

 \begin{table}[h!]
\centering
\caption{Characteristics of the images used for the image compression experiment.}\label{tab:images}
\begin{tabular}{lcc}
\hline
Image & Resolution & p \\
\hline
Atacama & 256$\times$192 & 100 \\
Santiago & 300$\times$250 & 100 \\
Valdivia & 300$\times$250 & 100 \\
Niebla & 400$\times$300& 150\\
\hline
\end{tabular}
\end{table}

{We shall implement ARBPG with blocks of size $b\in\{1,5,10,20\}$. Specifically, each coordinate block $x_i$, for $i=1,\ldots,2 \lceil p/b\rceil$ (with $\lceil \cdot \rceil$ denoting the ceiling function), compresses either $n_i = b \times m$ individual variables corresponding to $b$ rows of $V$, or $n_i = b \times n$ variables corresponding to $b$ rows of $W$. In particular, for the ``Niebla'' image with $b=20$, the last two coordinate blocks are determined by the last 10 rows of $V$ and $W$, respectively. In addition, we also run PALM and iPALM, the latter with inertial parameters given by $0.2$, which is the selection used in~\cite{MR3572363}.\footnote{In~\cite{MR3572363}, it is reported that a dynamical choice of the inertial parameters in the spirit of Nesterov's acceleration provides substantially better numerical results. This was also observed in our preliminary tests, but we do not include it here since no convergence guarantees are known.}   We choose the stopping criterion of the algorithms so that the same number of variable updates are measured.  Specifically, we stop ARBPG when the relative residual in a range of $2\lceil p/ b\rceil$ iterations is smaller than or equal to $10^{-4}$, i.e.,
\begin{equation}\label{eq:stopcrit}
\frac{\bigr| \varphi(\mU^{\hat{k}}, \mV^{\hat{k}}) - \varphi(\mU^k,\mV^k) \bigr|}{\|A\|_{F}} \leq 10^{-4}, 
\end{equation} 
where $\hat{k} = \max \{ k-2\lceil p/b\rceil,0\}$. For PALM and iPALM, we just measure the relative error between consecutive iterations, as these methods update all variables at every iteration. }

\begin{figure}[h!]
    \centering
    \includegraphics[width=0.98\textwidth]{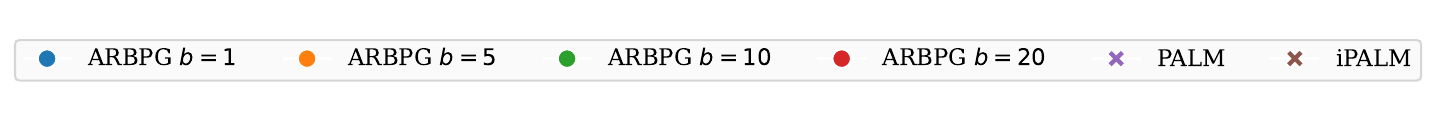} \\
    \includegraphics[width=0.98\textwidth]{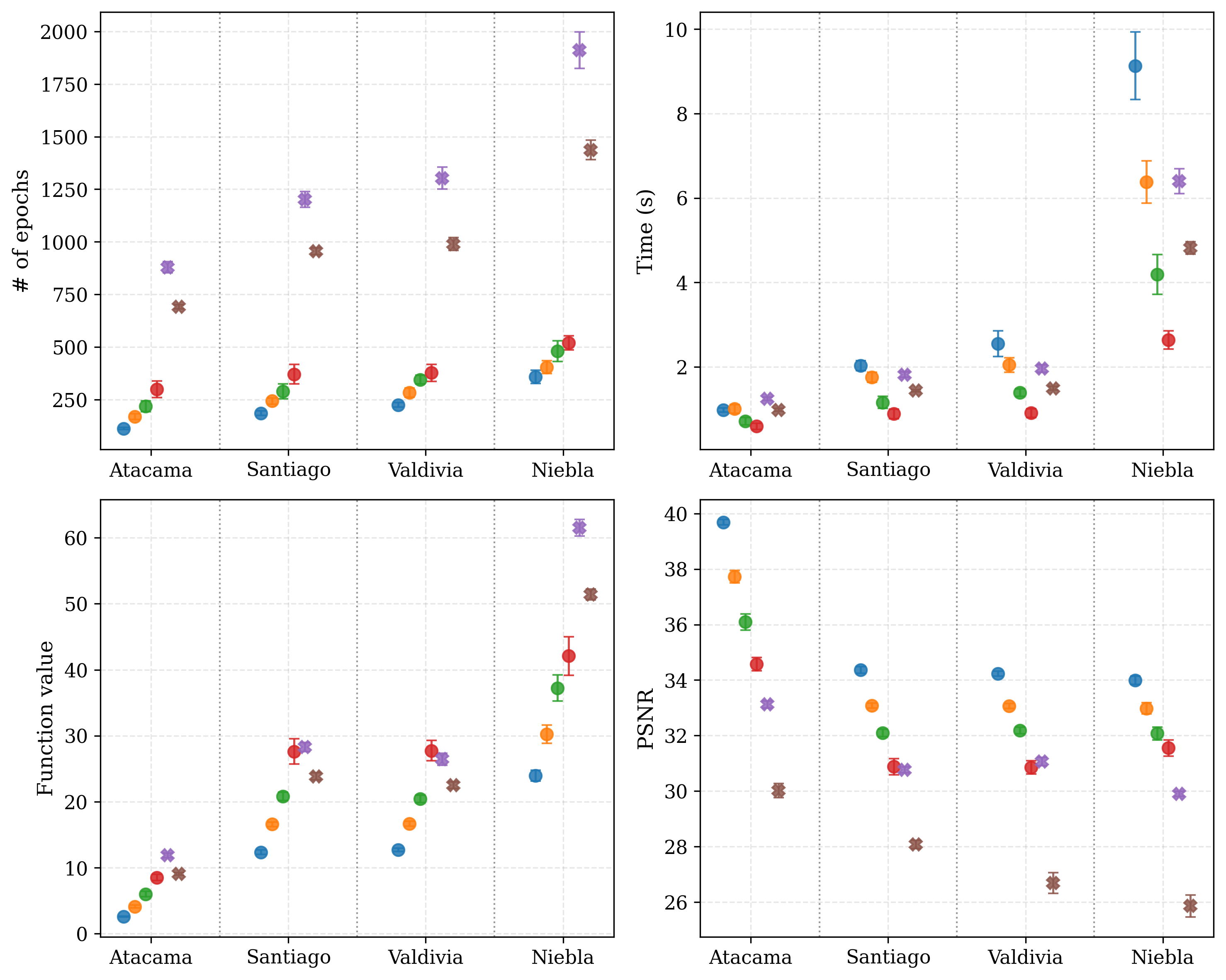}
    \caption{{Performance comparison of the evaluated algorithms across four test images (``Atacama", ``Santiago", ``Valdivia", and ``Niebla") over ten randomly generated seeds. The panels display the average number of epochs, computational time in seconds, final function values, and PSNR. The error bars denote the standard deviation across the ten independent runs.} }
\label{fig:compression1}
\end{figure}

{Our first test aims to compare the  performance of the different choices of block sizes for ARBPG. To this purpose, we consider the red color channel of each image, and run ARBPG, PALM and iPALM, initialized at the same ten pairs of  factor matrices $\mU\in\R^{p\times m}$ and $\mV\in\R^{p\times n}$  which are randomly generated, for each picture, with entries uniformly distributed in the interval $[0,1]$. As ARBPG and iPALM update a different number of variables at each iteration, this comparison is meaningless. Instead, we will compare the number of \emph{epochs}, that is, the number of  iterations required to update as many variables as present in the model.  For iPALM, an epoch coincides with an iteration of the algorithm, whereas for ARBPG an epoch is given by $2\lceil p/b\rceil $ iterations. Figure~\ref{fig:compression1} displays the average number of epochs, running time in seconds, and function values in the output. The bottom right plot   reports the average \emph{peak signal-to-noise ratio} (PSNR), which is a measure of the fidelity of the representation of the compressed image. The PSNR is measured in decibels and is given by
\[
PSNR = 20\, \mathrm{log}_{10}\left(MAX_A\right)  - 10\, \mathrm{log}_{10} \left( \frac{1}{m \cdot n}  \bigl\|A - (\mU^{out})^{T} \mV^{out}\bigr\|_F^2 \right),
\]
where $MAX_A$ is the maximum pixel value in $A$ and $\mU^{out}$, $\mV^{out}$ are the factor matrices returned by the algorithm at the stopping. The higher the value of the PSNR, the better the quality of compression. 
}

{Our results report a clear correlation between the block size $b$ and  
the number of epochs, the final function value, and the PSNR. Specifically, smaller 
block sizes lead to better performance across these three measures: fewer epochs, 
lower function values, and higher PSNR. The only exception is observed for ARBPG 
with $b=20$, which occasionally underperforms compared to PALM in terms of function 
value and PSNR. This relationship is inverted with respect to the computational cost for ARBPG, where 
larger values of $b$ significantly reduce the running time. In particular, PALM requires similar running time to ARBPG with $b=5$, but its 
PSNR is considerably worse. Furthermore, it is worth noting that in this experiment 
PALM consistently outperformed iPALM---under the chosen inertial configuration---in every 
evaluated aspect. On the whole, ARBPG with $b=1$ can be considered the best performing algorithm, 
as it yields superior PSNR values with only a slight increase in time (except for the ``Niebla'' image, where the time increase relative to $b=5$ is more pronounced).}

\begin{figure}[htbp]
    \centering
    \setlength{\arrayrulewidth}{0.5mm} 
    \setlength{\tabcolsep}{2pt} 
    \renewcommand{\arraystretch}{2} 

    \begin{tabular}{c c c c c} 
        & \textbf{Atacama} & \textbf{Santiago} & \textbf{Valdivia} & \textbf{Niebla} \\
   
   	\raisebox{28pt}{\rotatebox[origin=c]{90}{\textbf{Original}}} & 
        \includegraphics[width=.23\textwidth]{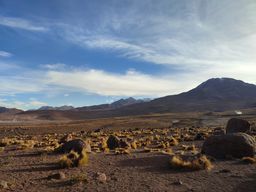} & 
        \includegraphics[width=.23\textwidth]{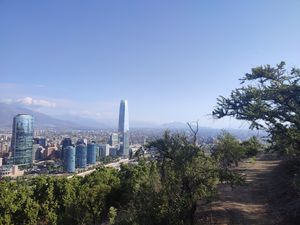} & 
        \includegraphics[width=.23\textwidth]{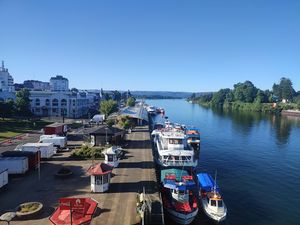} &
        \includegraphics[width=.23\textwidth]{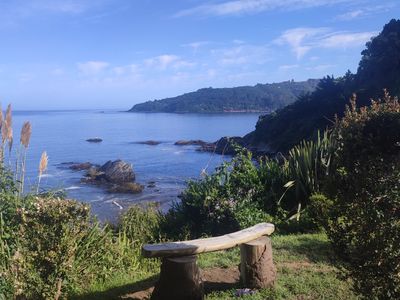} \\

        \raisebox{28pt}{\rotatebox[origin=c]{90}{\textbf{ARBPG $b=1$}}} & 
        \includegraphics[width=.23\textwidth]{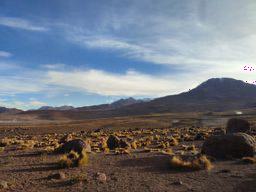} & 
        \includegraphics[width=.23\textwidth]{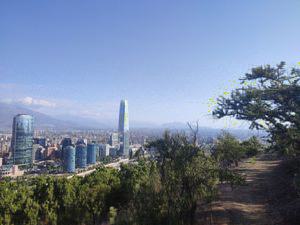} & 
        \includegraphics[width=.23\textwidth]{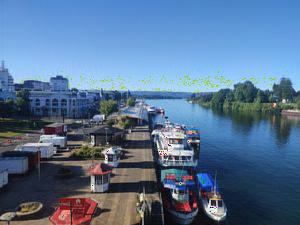} &
        \includegraphics[width=.23\textwidth]{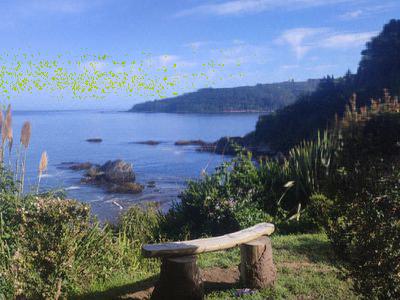}  \\

        \raisebox{28pt}{\rotatebox[origin=c]{90}{\textbf{ARBPG $b=5$}}} & 
        \includegraphics[width=.23\textwidth]{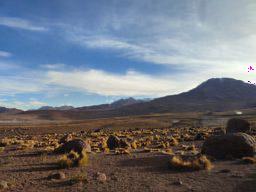} & 
        \includegraphics[width=.23\textwidth]{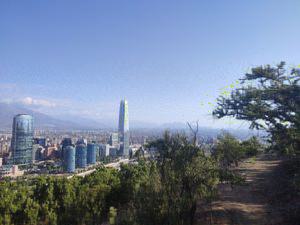} & 
        \includegraphics[width=.23\textwidth]{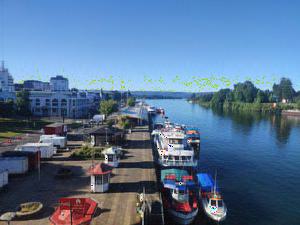} &
        \includegraphics[width=.23\textwidth]{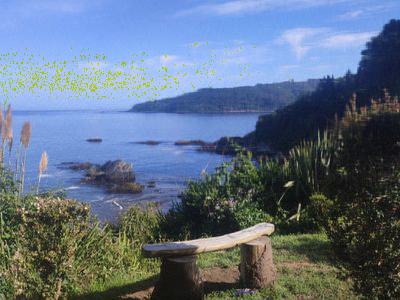} \\

        \raisebox{27pt}{\rotatebox[origin=c]{90}{\textbf{PALM}}} & 
        \includegraphics[width=.23\textwidth]{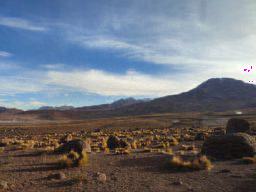} & 
        \includegraphics[width=.23\textwidth]{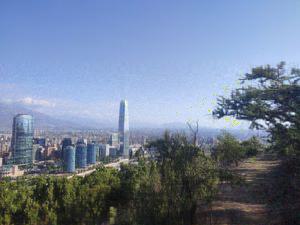} & 
        \includegraphics[width=.23\textwidth]{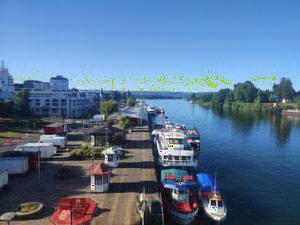} & 
        \includegraphics[width=.23\textwidth]{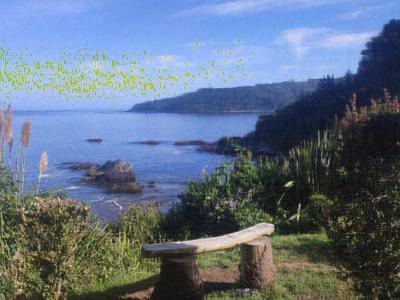}  \\

        \raisebox{27pt}{\rotatebox[origin=c]{90}{\textbf{SCIKIT}}} & 
        \includegraphics[width=.23\textwidth]{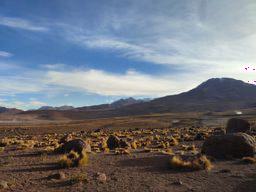} & 
        \includegraphics[width=.23\textwidth]{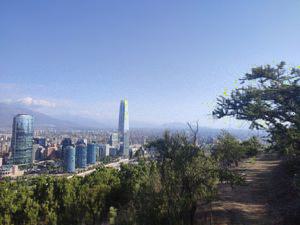} & 
        \includegraphics[width=.23\textwidth]{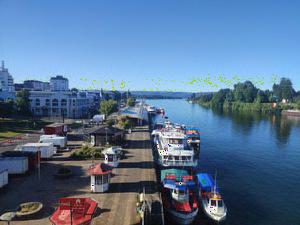} & 
        \includegraphics[width=.23\textwidth]{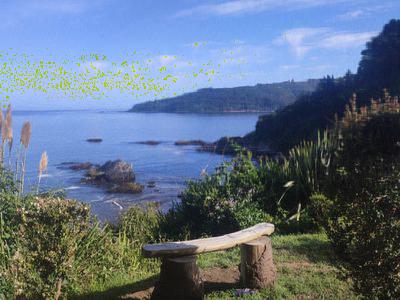}  \\
    \end{tabular}
    \caption{Original and  low rank compression of four Chilean landscapes images obtained by each method.}
\label{fig:image_compression}
\end{figure}

{We now present the results of carrying out the complete compression of the images. According to the results observed in Figure~\ref{fig:compression1}, we restrict our test to ARBPG with $b\in\{1,5\}$ and PALM, as these methods obtained the best PSNR values of their respective classes. In addition, we also consider \texttt{scikit-learn} Python library~\cite{scikit} built-in coordinate-descent NMF solver, which we will refer to as SCIKIT from now on. Since SCIKIT is very sensitive to initialization, we input into all  algorithms  the matrices returned by the warm start strategy of SCIKIT. On the other hand, to avoid excessive running time of SCIKIT, we set as its stopping criterion the same number of epochs as ARBPG $b=5$.
 Figure~\ref{fig:image_compression} illustrates the quality of compression obtained by each of the methods for the different images. In Table~\ref{t:images_results}, we show the obtained number of iterations, running time (in seconds), the function values at the output $\bigl(\varphi(V^{out},W^{out})\bigr)$ and the PSNR of each method, averaged over the three color channels (red, green, and blue) for every image.  The ARBPG methods require up to 5 to 7 times fewer epochs and similar running times than PALM, while consistently achieving 
lower function values and higher PSNR.   SCIKIT is undoubtedly the best performing method: it matches the PSNR of ARBPG with $b=1$, while requiring the least running time. Nonetheless, we emphasize that these results are heavily dependent on the 
initialization strategy employed for SCIKIT. Indeed, when initialized with random matrices, SCIKIT usually fails to achieve high-quality compression; 
in contrast, both ARBPG and PALM exhibit significantly greater robustness to the choice of the initial point.
}

\begin{table}[h!]

\centering
\caption{Summary of the results of the performance of each method applied for  the nonnegative matrix factorization  compression. The table displays the average value over the three color channels (red, green and blue).}
\label{t:images_results}
\begin{tblr}{width = \linewidth,colspec = {QQS[table-format=4.1,detect-weight,  mode=text]S[table-format=3.2, detect-weight,  mode=text]S[table-format=3.2, detect-weight,  mode=text]S[table-format=2.2, detect-weight,  mode=text]},}
\hline
{{{Image}}} & Method & {{{\# epochs}}}  & {{{Time (s)}}} & {{{$\varphi(V^{out},W^{out})$}}} & {{{PSNR}}} \\
\hline
Atacama& ARBPG $b=1$ &80 & 0.71 & 2.87&39.33\\
& ARBPG $b=5$ &111 & 0.63 & 4.02&37.88\\
& PALM &609 & 0.88 & 9.52&34.17\\
& SCIKIT &111 & 0.25 & 2.47&40.00\\
\hline Santiago & ARBPG $b=1$ &125 & 1.37 & 11.54&34.67\\
& ARBPG $b=5$ &151 & 1.15 & 15.69&33.33\\
& PALM &866 & 1.35 & 24.49&31.40\\
& SCIKIT &151 & 0.41 & 11.5&34.68\\
\hline Valdivia & ARBPG $b=1$ &153 & 1.69 & 12.18&34.43\\
& ARBPG $b=5$ &193 & 1.36 & 15.85&33.29\\
& PALM &969 & 1.49 & 24.85&31.34\\
& SCIKIT &193 & 0.51 & 11.81&34.56\\
\hline Niebla & ARBPG $b=1$ &277 & 7.03 & 22.96&34.10\\
& ARBPG $b=5$ &336 & 5.23 & 30.76&32.83\\
& PALM &1426 & 4.98 & 54.24&30.38\\
& SCIKIT &336 & 2.77 & 22.95&34.11\\
\hline
\end{tblr}
\end{table}

{
\subsubsection{Learning the  parts of objects }
The seminal work of Lee and Seung~\cite{LeeNMF,LeeSeung2000} demonstrated the ability of nonnegative matrix factorization to decompose nonnegative data in terms of parts-based representations. More specifically, assume the matrix $A\in\R^{m\times n}$ is decomposed into two nonnegative matrices $\mU\in\R^{p\times m}_+$ and $\mV\in\R^{p\times n}_+$, as $A\approx \mU^T \mV$. Then, the $p$ rows of $\mU$, which are the so-called \emph{basis vectors}, are expected to identify intrinsic ``parts" underlying  the original data.

In this section, we use our algorithm to illustrate this phenomenon. To the nonnegative factorization model considered in the previous section, we now incorporate an additional $\ell_0$-norm sparse constraint, that aims to enforce more compact basis representations. That is, we apply our algorithm to the mathematical program
\[
\min_{\mU,\mV} \;   \frac{1}{2} \|A-\mU^T\mV\|^2_F \quad \text{s.t.} \quad  \mU\in \R^{p\times m}_+;  \mV\in\R^{p\times n}_+; \|\mU_{i}\|_0 \leq s, i=1,\ldots,p;
\]
where $s$ is an integer and $\mU_{i}$ is the $i$-th row of $\mU$. Namely, the basis vectors $\mU_{i}$ should have at most $s$ nonzero elements. 
We now consider the $2$ blocks given by the complete factor matrix $\mU$ and $\mV$. Note that the row-wise projection of a vector $V$ onto the set $\{v\in\R^{m}_+,  \|v\|_0 \leq s\}$  can be efficiently computed with  closed formula given in~\cite[Proposition~4]{MR3232623}. Our database replicates the synthetic \emph{swimmer database} proposed by Donoho and Stodden in~\cite{DonohoStodden2003}.  Such a database consists of binary images of $20\times20$ pixels each that represent a swimmer composed of a motionless torso and the four limbs that can point in four different dimensions. The 256 resulting images are displayed in Figure~\ref{fig:collage}.  We ran ARBPG  to find $p=16$ basis images with sparsity less than $s=132\approx 0.33\times 400$. As observed in Figure~\ref{fig:basis_swimmer}, the basis images obtained by ARBPG successfully identify all possible positions of the limbs as ``parts-based building blocks" for constructing the entire collection of images. {Our numerical tests along multiple random initializations showed that the behavior of ARBPG for this  problem is very similar in the number of epochs, the running time and the quality of the solution to iPALM with $\alpha=\beta=0.2$, which in turn reported better results than the basic PALM algorithm.}


\begin{figure}[htbp]
    \centering
\includegraphics[width=\textwidth]{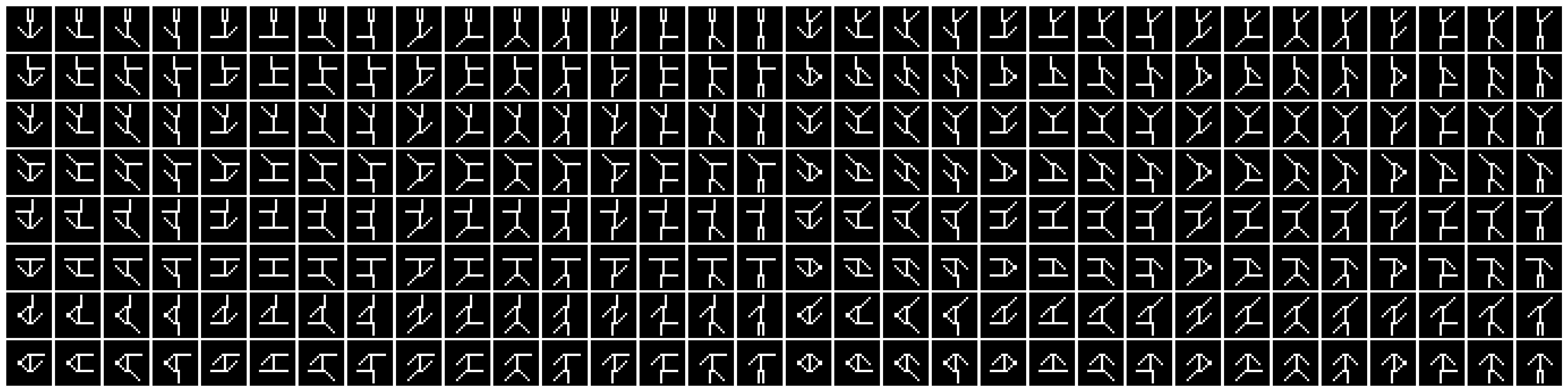}
\caption{The 256 pictures from the synthetic \emph{swimmer database}.}
\label{fig:collage}
\end{figure}



\begin{figure}
 \centering
 \includegraphics[width=.48\textwidth]{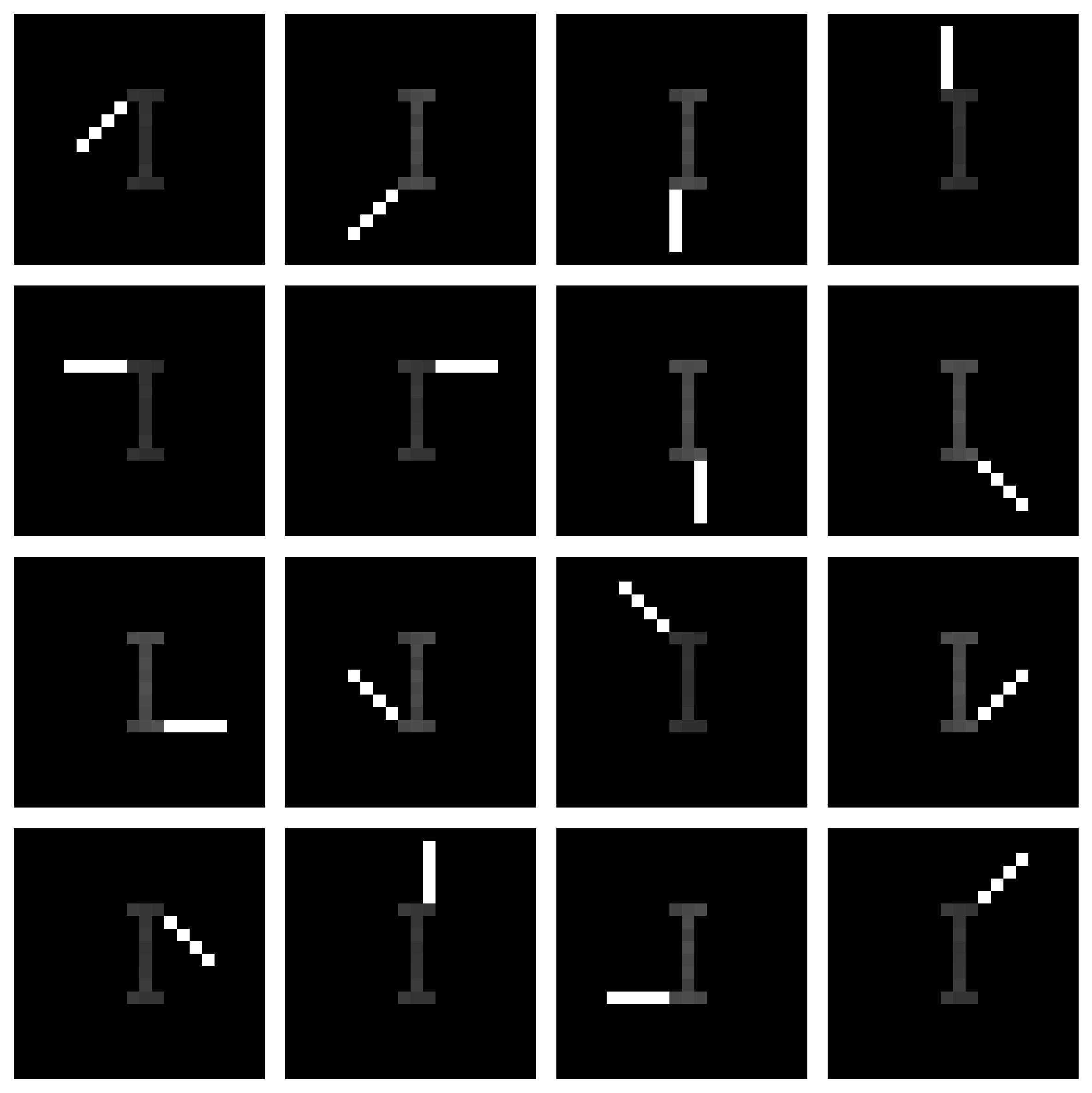}
 \caption{Basis images with sparsity $s=132$  of the swimmer dataset obtained after applying ARBPG.} \label{fig:basis_swimmer}
\end{figure}


{
\subsection{Symmetric nonnegative matrix factorization}\label{sect:num3}

Although NMF was originally introduced as a dimensionality reduction technique, it has proved successful as a clustering method in a wide variety of applications. Recently, several works have introduced a variant of the NMF model specifically tailored for complex graph and image clustering tasks, where the classical coordinate-space NMF approach is known to struggle; see, e.g.,~\cite{MR3353930,DingSymNMF,MR2785131} and the references therein. This formulation, known as \emph{Symmetric Nonnegative Matrix Factorization} (SymNMF), constructs a low-rank representation from a graph affinity matrix and is expressed as the optimization problem 
\begin{equation}\label{eq:SymNMF}
\min_{H} \quad \frac{1}{2}\| D - H^TH\|^2_F \quad \text{s.t.} \quad H\in\mathbb{R}^{p\times n}_+,
\end{equation}
where $D\in\mathbb{S}^n$ is a symmetric, nonnegative matrix whose entries represent pairwise similarity values between the $n$ data items, and $p$ denotes the target clustering rank.

From an optimization point of view, SymNMF exhibits significant differences with respect to the standard NMF problem~\eqref{eq:NMF1}. Let $f(H) = \frac{1}{2} \|D-H^TH\|^2_F$, and let $H_i \in \mathbb{R}^n$, for $i \in \{1,\ldots,p\}$, be the column vector representing the $i$-th row of~$H$. The (partial) gradients of $f$ with respect to $H$ and the block coordinate $H_i$ are given, respectively, by
\[
\nabla_H f(H) = -2 H (D-H^TH) \quad \text{and} \quad \nabla_{H_i} f(H) = - 2H_i^T (D-H^TH).
\]
Unlike the standard NMF gradients in~\eqref{eq:NMF_grad}, these partial gradients are no longer linear in the  variables of differentiation. This leads to two crucial distinctions with respect to our previous NMF experiment. First, standard PALM-type methods are not applicable here because, when fixed the non-differentiation variable, these partial gradients are only locally, rather than globally, Lipschitz continuous. Second, their local Lipschitz moduli cannot be efficiently computed exactly by evaluating the  norm of the second derivative---a rank-4 tensor whose computation should be avoided in practice---which is a desirable property when establishing an appropriate stepsize parameter (see Remark~\ref{remark:trial_step}). 

These challenges directly motivate the use of ARBPG with blocks defined by the individual rows of $H$. In this coordinate context, the second partial derivative of $f$ with respect to the row-block vector $H_i$ is given by the $n \times n$ matrix
\begin{equation*}
\nabla^2_{H_i} f(H) = -2(D-H^TH) + 2H_iH_i^T + 2\|H_i\|^2 I_n .
\end{equation*}
Constructing this full matrix and its spectral 2-norm at every single iteration to choose a trial stepsize is computationally expensive. Nonetheless, the triangle inequality provides an easily computable upper approximation of its spectral norm, given by
\begin{equation}\label{eq:boundsn}
\tilde{L}(H) := 2 \bigl( 2f(H) \bigr)^{1/2} + 4\|H_i\|^2,
\end{equation}
which exclusively requires cheap vector and scalar operations. As a consequence, a robust trial stepsize for ARBPG can be approximated as $\bar{\tau}_k = 1/\tilde{L}(H^k)$ for every $k\in\mathbb{N}$. Since $\tilde{L}(H^k)$ may provide a loose upper approximation of the true Lipschitz constant, this is the perfect setting to investigate whether the boosted linesearch provides an enhanced performance. 

In this spirit, we consider the \emph{digits dataset} contained in the \texttt{scikit-learn} Python library~\cite{scikit}. This data set contains $8\times 8$ pixel images of hand-written numbers between 0 and~9.
For constructing the similarity matrix,  we follow the recommendation in~\cite{MR2409803} and consider the $K$-nearest neighbor graph with $K=\lceil \log(n)\rceil +1$, where  $n$ is the number of data items. Specifically,  for each $n\in\{500, 800, 1000\}$ and 10 repetitions each, we randomly choose   $n$ images  from the dataset and initialize both 
ARBPG and ARBPG-B at the same matrix $H\in\mathbb{R}^{10\times n}$. In addition, the same seed is taken for both algorithms, so that the index updated at the same iteration coincides. We first run ARBPG using the same stopping condition as explained in~\eqref{eq:stopcrit}. Subsequently, we run ARBPG-B during the same execution  time used by ARBPG. 

The experimental results are reported in Table~\ref{t:SymNMF}. Columns 3 to 5   display the mean and standard deviation over the 10 instances for each dimension, of the number of iterations, number of function evaluations, and the number of boosted stepsizes accepted. The sixth column shows the \emph{clustering accuracy}, that is, the percentage of data points correctly classified by the algorithm relative to the ground-truth labels. The seventh column gathers the \emph{adjusted rand index} (ARI) that measures the quality of the clustering by counting pairs that are assigned in the same or different clusters in the predicted and true clusterings. Perfect labeling achieves an ARI of 1, whereas poorer agreeing assignments yield lower scores up to $-0.5$, with  0 corresponding to random labeling.

\begin{table}[h!]
\centering
\caption{Performance comparison between ARBPG and ARBPG-B on the handwritten digits dataset across varying sample dimensions ($n \in \{500, 800, 1000\}$). Results report the mean $\pm$ standard deviation computed over 10 independent randomized trials. Performance is evaluated across the number of iterations, objective function, and number of boosted stepsizes accepted by ARBPG-B, alongside two clustering metrics, namely,  clustering accuracy and adjusted rand index.}
\label{t:SymNMF}
\begin{tblr}{width = \linewidth, colsep = 2.5pt,  rows = {abovesep=3pt, belowsep=3pt}, colspec = {Q[c, valign=m] Q[l, valign=m] X[c, valign=m] X[c, valign=m] X[c, valign=m] X[c, valign=m] X[c, valign=m]},}
\hline
Dim. & Method & {\# iterations} & {\# eval. $\varphi$} & {\# boosts} & {Clust. acc. (\%)} & {ARI} \\
\hline
   500 & ARBPG & $3483 \pm 737$ & $3483 \pm 737$ &  & $75.66 \pm 6.59$ & $0.66 \pm 0.07$ \\
 & ARBPG-B & $2768 \pm 511$ & $5538 \pm 1024$ & $1043 \pm 367$ & $77.30 \pm 6.57$ & $0.69 \pm 0.08$ \\
\hline
  800 & ARBPG & $5066 \pm 1374$ & $5066 \pm 1374$ &  & $68.50 \pm 10.50$ & $0.60 \pm 0.11$ \\
 & ARBPG-B & $4115 \pm 1016$ & $8226 \pm 2022$ & $1530 \pm 607$ & $71.72 \pm 10.08$ & $0.62 \pm 0.11$ \\
\hline
 1000 & ARBPG & $4394 \pm 559$ & $4394 \pm 559$ &  & $66.42 \pm 7.25$ & $0.56 \pm 0.06$ \\
 & ARBPG-B & $3701 \pm 415$ & $7403 \pm 829$ & $1077 \pm 248$ & $68.78 \pm 7.90$ & $0.58 \pm 0.08$ \\
\hline
\end{tblr}
\end{table}

A primary observation is that ARBPG requires the same number of iterations and function evaluations. This confirms that \eqref{eq:boundsn}  is a valid upper approximation of the local Lipschitz modulus. Nonetheless, the approximation is not tight, which leads to the boosted linesearch being accepted in approximately 30\% of the iterations. As a consequence, ARBPG-B consistently and significantly reduces the total number of main iterations required for convergence across all tested scales. We report the average and maximum boosted stepsizes $\lambda_k$ (including the instances when $\lambda_k=0$) accepted by ARBPG-B in Table~\ref{t:SymNMF_boost}. Naturally, the activation of the boosted linesearch leads to an increase in the number of function evaluations relative to ARBPG. However, given that these extra calculations exclusively involve cheap scalar and vector operations and avoid forming full-size matrices, the overall computational overhead per evaluation remains low. The running times of the algorithms were on average 10.16 seconds for $n=500$, 39.86 seconds for $n=800$ and 58.45 seconds for $n=1000$.

\begin{table}[h!]

\centering
\caption{For the ARBPG-B runs in Table~\ref{t:SymNMF},  average and maximum values of the boosted stepsize $\lambda_k$ (including 0 values) averaged across the 10 trials for each  sample dimension $n \in \{500, 800, 1000\}$.}
\label{t:SymNMF_boost}
\begin{tblr}{width = .5\linewidth, colsep = 2.5pt,  rows = {abovesep=3pt, belowsep=3pt}, colspec = {Q[c, valign=m]  X[c, valign=m] X[c, valign=m]},}
\hline
Dim.  & {average $\lambda_k$} & {max $\lambda_k$}  \\
\hline
500  & $1.00 \pm 0.42$ & $37.40 \pm 17.87$  \\
\hline
800  & $0.92 \pm 0.37$ & $45.40 \pm 33.14$  \\
\hline
1000  & $0.54 \pm 0.13$ & $26.20 \pm 14.40$  \\
\hline
\end{tblr}
\end{table}

Regarding the quality of the clustering, ARBPG-B exhibits a consistent improvement in  performance for all dimensions. At $n=500$, the boosted variant achieves a higher clustering accuracy ($77.30\% \pm 6.57\%$) and ARI ($0.69 \pm 0.08$) compared to standard ARBPG ($75.66\% \pm 6.59\%$ accuracy and $0.66 \pm 0.07$ ARI). Similarly, for $n=800$ and $n=1000$, ARBPG-B yields an improvement of more than $2$ percentage points in accuracy and $0.02$ in ARI. In particular, only in 3 out of the 30 trials ARBPG obtained a higher clustering accuracy than ARBPG-B. Furthermore, for $n=1000$ ARBPG-B obtained better clustering performance in every instance.  On the whole, our results showcase that the use of ARBPG-B  for this problem is advantageous, as  the same algorithmic execution time provides noticeable improvement in clustering results.}

\section{Conclusions and future work}\label{sect:conclusion}

{This paper extended
the framework of randomized block-coordinate descent algorithms to mathematical programs lacking the usual assumption of global Lipschitz continuity of  block partial gradients of smooth functions.} Specifically, we considered the nonconvex optimization problem of minimizing a separable lower semicontinuous function and a smooth function whose gradient is block-wise  Lipschitz continuous only locally. We proposed an Adaptive Randomized Block Proximal Gradient (ARBPG) method that adaptively selects the proximal stepsize to ensure a sufficient decrease condition without requiring knowledge of the local Lipschitz moduli of the partial  gradients of the smooth function. In addition, ARBPG allows the possibility of performing a boosted linesearch in the spirit of~\cite{BDSA2025} with a self-adaptive procedure to update the parameters. {Our numerical experiments show that ARBPG competes well with the widely-known PALM and iPALM methods for the Nonnegative Matrix Factorization (NMF) problem. Furthermore, the boosted version of ARBPG, called ARBPG-B, was reported to be advantageous in the symmetric NMF problem, where inexact estimations of the local Lipschitz moduli led to poor choices of proximal stepsizes.}

In a future work, we would like to investigate the convergence of Algorithm~\ref{alg:1} under some error bound conditions.  For instance, in~\cite{MR3293507} a local error bound condition was used to establish a linear convergence rate for the sequence of function values generated by a randomized block proximal gradient algorithm. However, they needed to impose the stronger assumption of global Lipschitz continuity of the whole gradient. The extension of these techniques to cover the local Lipschitz gradient case and the adaptive procedure conducted by Algorithm~\ref{alg:1}  seems a challenging task.

\paragraph{Data availability} The Python source code and data of the
experiments are available at \href{https://github.com/DavidTBelen/ARBPG}{https://github.com/DavidTBelen/ARBPG}.

\section*{Statements and declarations}

\paragraph{Acknowledgment} 
The authors are grateful to the anonymous referees for their valuable comments and suggestions, which significantly improved the quality of this paper.

\paragraph{Competing interests} 
The authors declare that they have no competing interests.

\bibliographystyle{acm}
\bibliography{references}

@article {wright2020analyzing,
    AUTHOR = {Wright, Stephen J. and Lee, Ching-pei},
     TITLE = {Analyzing random permutations for cyclic coordinate descent},
   JOURNAL = {Math. Comp.},
  FJOURNAL = {Mathematics of Computation},
    VOLUME = {89},
      YEAR = {2020},
    NUMBER = {325},
     PAGES = {2217--2248},
      ISSN = {0025-5718},
   MRCLASS = {65K05 (68W20 90C25)},
  MRNUMBER = {4109565},
MRREVIEWER = {Hideaki Iiduka},
       DOI = {10.1090/mcom/3530},
       URL = {https://doi.org/10.1090/mcom/3530},
}

@article {gurbuzbalaban2020randomness,
    AUTHOR = {G\"{u}rb\"{u}zbalaban, Mert and Ozdaglar, Asuman and Vanli, Nuri
              Denizcan and Wright, Stephen J.},
     TITLE = {Randomness and permutations in coordinate descent methods},
   JOURNAL = {Math. Program.},
  FJOURNAL = {Mathematical Programming},
    VOLUME = {181},
      YEAR = {2020},
    NUMBER = {2, Ser. B},
     PAGES = {349--376},
      ISSN = {0025-5610},
   MRCLASS = {90C06 (65K05 90C25)},
  MRNUMBER = {4109506},
MRREVIEWER = {Miantao Chao},
       DOI = {10.1007/s10107-019-01438-4},
       URL = {https://doi.org/10.1007/s10107-019-01438-4},
}

@book {MR3719240,
    AUTHOR = {Beck, Amir},
     TITLE = {First-order methods in optimization},
    SERIES = {MOS-SIAM Series on Optimization},
    VOLUME = {25},
 PUBLISHER = {Society for Industrial and Applied Mathematics (SIAM),
              Philadelphia, PA; Mathematical Optimization Society,
              Philadelphia, PA},
      YEAR = {2017},
     PAGES = {xii+475},
      ISBN = {978-1-611974-98-0},
   MRCLASS = {90-02 (49M20 65K05 90C06 90C25)},
  MRNUMBER = {3719240},
MRREVIEWER = {Sven-\AA ke\ Gustafson},
       DOI = {10.1137/1.9781611974997.ch1},
       URL = {https://doi.org/10.1137/1.9781611974997.ch1},
}

@book {MR4704094,
	AUTHOR = {Tao, Terence},
	TITLE = {Analysis {II}},
	SERIES = {Texts and Readings in Mathematics},
	VOLUME = {38},
	EDITION = {Fourth},
	PUBLISHER = {Springer, Singapore; Hindustan Book Agency, New Delhi},
	YEAR = {[2022] \copyright 2022},
	PAGES = {xvii+195},
	ISBN = {978-9-81197-284-3},
	MRCLASS = {26-01},
	MRNUMBER = {4704094},
	DOI = {10.1007/978-981-19-7284-3},
	URL = {https://doi.org/10.1007/978-981-19-7284-3},
}

@book {MR270403,
	AUTHOR = {Feller, William},
	TITLE = {An introduction to probability theory and its applications.
	{V}ol. {II}},
	EDITION = {Second},
	PUBLISHER = {John Wiley \& Sons, Inc., New York-London-Sydney},
	YEAR = {1971},
	PAGES = {xxiv+669},
	MRCLASS = {60.00},
	MRNUMBER = {270403},
}

@book {MR2722836,
	AUTHOR = {Durrett, Rick},
	TITLE = {Probability: theory and examples},
	SERIES = {Cambridge Series in Statistical and Probabilistic Mathematics},
	VOLUME = {31},
	EDITION = {Fourth},
	PUBLISHER = {Cambridge University Press, Cambridge},
	YEAR = {2010},
	PAGES = {x+428},
	ISBN = {978-0-521-76539-8},
	MRCLASS = {60-01},
	MRNUMBER = {2722836},
	DOI = {10.1017/CBO9780511779398},
	URL = {https://doi.org/10.1017/CBO9780511779398},
}

@book {MR2267655,
	AUTHOR = {Bogachev, V. I.},
	TITLE = {Measure theory. {V}ol. {I}, {II}},
	PUBLISHER = {Springer-Verlag, Berlin},
	YEAR = {2007},
	PAGES = {Vol. I: xviii+500 pp., Vol. II: xiv+575},
	ISBN = {978-3-540-34513-8; 3-540-34513-2},
	MRCLASS = {28-02 (28Axx 28Cxx 46G12 60G42 60G44)},
	MRNUMBER = {2267655},
	MRREVIEWER = {Ren\'e\ L.\ Schilling},
	DOI = {10.1007/978-3-540-34514-5},
	URL = {https://doi.org/10.1007/978-3-540-34514-5},
}

@article {MR2968857,
    AUTHOR = {Nesterov, Yu.},
     TITLE = {Efficiency of coordinate descent methods on huge-scale
              optimization problems},
   JOURNAL = {SIAM J. Optim.},
  FJOURNAL = {SIAM Journal on Optimization},
    VOLUME = {22},
      YEAR = {2012},
    NUMBER = {2},
     PAGES = {341--362},
      ISSN = {1052-6234},
   MRCLASS = {90C06 (90C25 90C47)},
  MRNUMBER = {2968857},
MRREVIEWER = {Jan-J. R\"{u}ckmann},
       DOI = {10.1137/100802001},
       URL = {https://doi.org/10.1137/100802001},
}

@article {MR3404687,
    AUTHOR = {Fercoq, Olivier and Richt\'{a}rik, Peter},
     TITLE = {Accelerated, parallel, and proximal coordinate descent},
   JOURNAL = {SIAM J. Optim.},
  FJOURNAL = {SIAM Journal on Optimization},
    VOLUME = {25},
      YEAR = {2015},
    NUMBER = {4},
     PAGES = {1997--2023},
      ISSN = {1052-6234},
   MRCLASS = {65K05 (49M27 65Y20 68Q25 68W10 68W20 90C25)},
  MRNUMBER = {3404687},
MRREVIEWER = {Aurea Mart\'{\i}nez},
       DOI = {10.1137/130949993},
       URL = {https://doi.org/10.1137/130949993},
}

@article {MR3293507,
    AUTHOR = {Patrascu, Andrei and Necoara, Ion},
     TITLE = {Efficient random coordinate descent algorithms for large-scale
              structured nonconvex optimization},
   JOURNAL = {J. Global Optim.},
  FJOURNAL = {Journal of Global Optimization. An International Journal
              Dealing with Theoretical and Computational Aspects of Seeking
              Global Optima and Their Applications in Science, Management
              and Engineering},
    VOLUME = {61},
      YEAR = {2015},
    NUMBER = {1},
     PAGES = {19--46},
      ISSN = {0925-5001},
   MRCLASS = {90C26 (90C56)},
  MRNUMBER = {3293507},
MRREVIEWER = {Julien Ugon},
       DOI = {10.1007/s10898-014-0151-9},
       URL = {https://doi.org/10.1007/s10898-014-0151-9},
}

@book {MR1491362,
	AUTHOR = {Rockafellar, R. T. and Wets, R. J.-B.},
	TITLE = {Variational Analysis},
	PUBLISHER = {Springer-Verlag, Berlin},
	YEAR = {1998},
	PAGES = {xiv+733},
	ISBN = {3-540-62772-3},
	MRCLASS = {49-02 (46N10 47N10 49J52 49K40 90C30)},
	MRNUMBER = {1491362},
	MRREVIEWER = {Francis H. Clarke},
}

@article {MR3730541,
    AUTHOR = {Lu, Zhaosong and Xiao, Lin},
     TITLE = {A randomized nonmonotone block proximal gradient method for a
              class of structured nonlinear programming},
   JOURNAL = {SIAM J. Numer. Anal.},
  FJOURNAL = {SIAM Journal on Numerical Analysis},
    VOLUME = {55},
      YEAR = {2017},
    NUMBER = {6},
     PAGES = {2930--2955},
      ISSN = {0036-1429},
   MRCLASS = {65K05 (90C06 90C26)},
  MRNUMBER = {3730541},
MRREVIEWER = {Wah June Leong},
       DOI = {10.1137/16M1110182},
       URL = {https://doi.org/10.1137/16M1110182},
}

@Book{bauschke2017,
  author    = {Bauschke, Heinz H. and Combettes, Patrick L.},
  publisher = {2nd edn., Springer, Cham},
  title     = {Convex Analysis and Monotone Operator Theory in Hilbert Spaces},
  year      = {2017},
}

@book {MR3823783,
	AUTHOR = {Mordukhovich, B. S.},
	TITLE = {Variational Analysis and Applications},
	PUBLISHER = {Springer, Cham},
	YEAR = {2018},
	PAGES = {xix+622},
	ISBN = {978-3-319-92773-2; 978-3-319-92775-6},
	MRCLASS = {49-02 (49J52 49J53 49K40 90-02 90C29 90C30 90C34)},
	MRNUMBER = {3823783},
	MRREVIEWER = {Amos Uderzo},
	DOI = {10.1007/978-3-319-92775-6},
	URL = {https://doi.org/10.1007/978-3-319-92775-6},
}

@book {MR2191744,
	AUTHOR = {Mordukhovich, Boris S.},
	TITLE = {Variational Analysis and Generalized Differentiation {I}},
	PUBLISHER = {Springer-Verlag, Berlin},
	YEAR = {2006},
	PAGES = {xxii+579},
	ISBN = {978-3-540-25437-9; 3-540-25437-4},
	MRCLASS = {49-02 (49J52 49J53 49K27 49K40 90C31)},
	MRNUMBER = {2191744},
	MRREVIEWER = {Lionel Thibault},
}

@article {MR3160324,
    AUTHOR = {Shechtman, Yoav and Beck, Amir and Eldar, Yonina C.},
     TITLE = {G{ESPAR}: {E}fficient phase retrieval of sparse signals},
   JOURNAL = {IEEE Trans. Signal Process.},
  FJOURNAL = {IEEE Transactions on Signal Processing},
    VOLUME = {62},
      YEAR = {2014},
    NUMBER = {4},
     PAGES = {928--938},
      ISSN = {1053-587X},
   MRCLASS = {94A12},
  MRNUMBER = {3160324},
MRREVIEWER = {K. B. Datta},
       DOI = {10.1109/TSP.2013.2297687},
       URL = {https://doi.org/10.1109/TSP.2013.2297687},
}

@article {MR3080197,
    AUTHOR = {Beck, Amir and Eldar, Yonina C.},
     TITLE = {Sparsity constrained nonlinear optimization: optimality
              conditions and algorithms},
   JOURNAL = {SIAM J. Optim.},
  FJOURNAL = {SIAM Journal on Optimization},
    VOLUME = {23},
      YEAR = {2013},
    NUMBER = {3},
     PAGES = {1480--1509},
      ISSN = {1052-6234},
   MRCLASS = {90C26 (90C46)},
  MRNUMBER = {3080197},
MRREVIEWER = {Magdalena Nockowska-Rosiak},
       DOI = {10.1137/120869778},
       URL = {https://doi.org/10.1137/120869778},
}

@article {MR3832977,
    AUTHOR = {Bolte, J\'{e}r\^{o}me and Sabach, Shoham and Teboulle, Marc and
              Vaisbourd, Yakov},
     TITLE = {First order methods beyond convexity and {L}ipschitz gradient
              continuity with applications to quadratic inverse problems},
   JOURNAL = {SIAM J. Optim.},
  FJOURNAL = {SIAM Journal on Optimization},
    VOLUME = {28},
      YEAR = {2018},
    NUMBER = {3},
     PAGES = {2131--2151},
      ISSN = {1052-6234},
   MRCLASS = {90C25 (26B25 49M27 52A41 65K05)},
  MRNUMBER = {3832977},
MRREVIEWER = {Nguyen Hieu Thao},
       DOI = {10.1137/17M1138558},
       URL = {https://doi.org/10.1137/17M1138558},
}

@article {LeeNMF,
    AUTHOR = {Lee, Daniel and Seung, H.},
     TITLE = {Learning parts of objects by non-negative matrix factorization},
   JOURNAL = {Nature},
  FJOURNAL = {Nature},
    VOLUME = {401},
      YEAR = {1999},
     PAGES = {788--791},
       DOI = {10.1038/44565},
       URL = {https://doi.org/10.1038/44565},
}

@inproceedings {MR2388467,
    AUTHOR = {Pauca, V. Paul and Shahnaz, Farial and Berry, Michael W. and
              Plemmons, Robert J.},
     TITLE = {Text mining using non-negative matrix factorizations},
 BOOKTITLE = {Proceedings of the {F}ourth {SIAM} {I}nternational
              {C}onference on {D}ata {M}ining},
     PAGES = {452--456},
 PUBLISHER = {SIAM, Philadelphia, PA},
      YEAR = {2004},
   MRCLASS = {Expansion},
  MRNUMBER = {2388467},
}

@article{PAUCA200629,
title = {Nonnegative matrix factorization for spectral data analysis},
journal = {Linear Algebra and its Applications},
volume = {416},
number = {1},
pages = {29-47},
year = {2006},
note = {Special Issue devoted to the Haifa 2005 conference on matrix theory},
issn = {0024-3795},
doi = {https://doi.org/10.1016/j.laa.2005.06.025},
url = {https://www.sciencedirect.com/science/article/pii/S002437950500340X},
author = {V. Paul Pauca and J. Piper and Robert J. Plemmons}
}

@article {MR3158055,
    AUTHOR = {Kim, Jingu and He, Yunlong and Park, Haesun},
     TITLE = {Algorithms for nonnegative matrix and tensor factorizations: a
              unified view based on block coordinate descent framework},
   JOURNAL = {J. Global Optim.},
  FJOURNAL = {Journal of Global Optimization. An International Journal
              Dealing with Theoretical and Computational Aspects of Seeking
              Global Optima and Their Applications in Science, Management
              and Engineering},
    VOLUME = {58},
      YEAR = {2014},
    NUMBER = {2},
     PAGES = {285--319},
      ISSN = {0925-5001},
   MRCLASS = {65F05 (15A23 90C90)},
  MRNUMBER = {3158055},
MRREVIEWER = {Bogdan Dumitrescu},
       DOI = {10.1007/s10898-013-0035-4},
       URL = {https://doi.org/10.1007/s10898-013-0035-4},
}

@article {BDSA2025,
    AUTHOR = {Arag\'{o}n-Artacho, Francisco J. and P\'{e}rez-Aros, Pedro and
              Torregrosa-Bel\'{e}n, David},
     TITLE = {The boosted double-proximal subgradient algorithm for
              nonconvex optimization},
   JOURNAL = {Math. Program.},
  FJOURNAL = {Mathematical Programming},
    VOLUME = {214},
      YEAR = {2025},
    NUMBER = {1-2, Ser. A},
     PAGES = {491--537},
      ISSN = {0025-5610,1436-4646},
   MRCLASS = {49J53 (65K05 90C26 90C30 90C46)},
  MRNUMBER = {4996113},
       DOI = {10.1007/s10107-024-02190-0},
       URL = {https://doi.org/10.1007/s10107-024-02190-0},
}

@article {MR4078808,
    AUTHOR = {Arag\'{o}n Artacho, Francisco J. and Vuong, Phan T.},
     TITLE = {The boosted difference of convex functions algorithm for
              nonsmooth functions},
   JOURNAL = {SIAM J. Optim.},
  FJOURNAL = {SIAM Journal on Optimization},
    VOLUME = {30},
      YEAR = {2020},
    NUMBER = {1},
     PAGES = {980--1006},
      ISSN = {1052-6234},
   MRCLASS = {65K05 (47N10 65K10 90C26)},
  MRNUMBER = {4078808},
MRREVIEWER = {Majid E. Abbasov},
       DOI = {10.1137/18M123339X},
       URL = {https://doi.org/10.1137/18M123339X},
}

@book {MR1099605,
    AUTHOR = {Polyak, Boris T.},
     TITLE = {Introduction to optimization},
    SERIES = {Translations Series in Mathematics and Engineering},
      NOTE = {Translated from the Russian,
              With a foreword by Dimitri P. Bertsekas},
 PUBLISHER = {Optimization Software, Inc., Publications Division, New York},
      YEAR = {1987},
     PAGES = {xxvii+438},
      ISBN = {0-911575-14-6},
   MRCLASS = {49-01 (65Kxx 90Cxx)},
  MRNUMBER = {1099605},
}

@article {MR3785672,
    AUTHOR = {Arag\'{o}n Artacho, Francisco J. and Fleming, Ronan M. T. and
              Vuong, Phan T.},
     TITLE = {Accelerating the {DC} algorithm for smooth functions},
   JOURNAL = {Math. Program.},
  FJOURNAL = {Mathematical Programming},
    VOLUME = {169},
      YEAR = {2018},
    NUMBER = {1, Ser. B},
     PAGES = {95--118},
      ISSN = {0025-5610},
   MRCLASS = {65K05 (65K10 90C26 92C42)},
  MRNUMBER = {3785672},
MRREVIEWER = {Dariusz Uci\'{n}ski},
       DOI = {10.1007/s10107-017-1180-1},
       URL = {https://doi.org/10.1007/s10107-017-1180-1},
}

@incollection {MR874369,
    AUTHOR = {Pham Dinh Tao and Souad, El Bernoussi},
     TITLE = {Algorithms for solving a class of nonconvex optimization
              problems. {M}ethods of subgradients},
 BOOKTITLE = {F{ERMAT} days 85: mathematics for optimization ({T}oulouse,
              1985)},
    SERIES = {North-Holland Math. Stud.},
    VOLUME = {129},
     PAGES = {249--271},
 PUBLISHER = {North-Holland, Amsterdam},
      YEAR = {1986},
   MRCLASS = {90C33},
  MRNUMBER = {874369},
MRREVIEWER = {Wenci Yu},
       DOI = {10.1016/S0304-0208(08)72402-2},
       URL = {https://doi.org/10.1016/S0304-0208(08)72402-2},
}

@article {MR4757561,
    AUTHOR = {Ferreira, O. P. and Santos, E. M. and Souza, J. C. O.},
     TITLE = {A boosted {DC} algorithm for non-differentiable {DC}
              components with non-monotone line search},
   JOURNAL = {Comput. Optim. Appl.},
  FJOURNAL = {Computational Optimization and Applications. An International
              Journal},
    VOLUME = {88},
      YEAR = {2024},
    NUMBER = {3},
     PAGES = {783--818},
      ISSN = {0926-6003,1573-2894},
   MRCLASS = {90C26 (47N10 65K05 65K10)},
  MRNUMBER = {4757561},
       DOI = {10.1007/s10589-024-00578-4},
       URL = {https://doi.org/10.1007/s10589-024-00578-4},
}

@article {MR4527568,
    AUTHOR = {Arag\'{o}n-Artacho, F. J. and Campoy, R. and Vuong, P. T.},
     TITLE = {The boosted {DC} algorithm for linearly constrained {DC}
              programming},
   JOURNAL = {Set-Valued Var. Anal.},
  FJOURNAL = {Set-Valued and Variational Analysis. Theory and Applications},
    VOLUME = {30},
      YEAR = {2022},
    NUMBER = {4},
     PAGES = {1265--1289},
      ISSN = {1877-0533},
   MRCLASS = {65K10 (65K05 90C20 90C26)},
  MRNUMBER = {4527568},
       DOI = {10.1007/s11228-022-00656-x},
       URL = {https://doi.org/10.1007/s11228-022-00656-x},
}

@Article{tran2023boosteddcalgorithmclustering,
      title={The Boosted {DC} Algorithm for Clustering with Constraints}, 
      author={Tuyen Tran and Kate Figenschou and Phan Tu Vuong},
      journal = {Preprint, arXiv:\href{https://arxiv.org/abs/2310.14148}{2310.14148} [math.OC], 2023},
}

@article {MR359795,
    AUTHOR = {Auslender, A.},
     TITLE = {M\'{e}thodes num\'{e}riques pour la d\'{e}composition et la minimisation
              de fonctions non diff\'{e}rentiables},
   JOURNAL = {Numer. Math.},
  FJOURNAL = {Numerische Mathematik},
    VOLUME = {18},
      YEAR = {1971/72},
     PAGES = {213--223},
      ISSN = {0029-599X},
   MRCLASS = {90C25},
  MRNUMBER = {359795},
MRREVIEWER = {P. J. Laurent},
       DOI = {10.1007/BF01397082},
       URL = {https://doi.org/10.1007/BF01397082},
}

@book {MR273810,
    AUTHOR = {Ortega, J. M. and Rheinboldt, W. C.},
     TITLE = {Iterative solution of nonlinear equations in several
              variables},
 PUBLISHER = {Academic Press, New York-London},
      YEAR = {1970},
     PAGES = {xx+572},
   MRCLASS = {65.50},
  MRNUMBER = {273810},
MRREVIEWER = {J. F. Traub},
}

@book {MR3444832,
    AUTHOR = {Bertsekas, Dimitri P.},
     TITLE = {Nonlinear programming},
    SERIES = {Athena Scientific Optimization and Computation Series},
   EDITION = {Second},
 PUBLISHER = {Athena Scientific, Belmont, MA},
      YEAR = {1999},
     PAGES = {xiv+777},
      ISBN = {1-886529-00-0},
   MRCLASS = {49-01 (90-01)},
  MRNUMBER = {3444832},
}

@ARTICLE{491321,
  author={Bouman, C.A. and Sauer, K.},
  journal={IEEE Transactions on Image Processing}, 
  title={A unified approach to statistical tomography using coordinate descent optimization}, 
  year={1996},
  volume={5},
  number={3},
  pages={480-492},
  doi={10.1109/83.491321}}

@article {MR2421312,
    AUTHOR = {Tseng, Paul and Yun, Sangwoon},
     TITLE = {A coordinate gradient descent method for nonsmooth separable
              minimization},
   JOURNAL = {Math. Program.},
  FJOURNAL = {Mathematical Programming. A Publication of the Mathematical
              Programming Society},
    VOLUME = {117},
      YEAR = {2009},
    NUMBER = {1-2, Ser. B},
     PAGES = {387--423},
      ISSN = {0025-5610},
   MRCLASS = {49M27 (65K10 90C06)},
  MRNUMBER = {2421312},
MRREVIEWER = {Marcos Raydan},
       DOI = {10.1007/s10107-007-0170-0},
       URL = {https://doi.org/10.1007/s10107-007-0170-0},
}

@article {MR3347548,
    AUTHOR = {Wright, Stephen J.},
     TITLE = {Coordinate descent algorithms},
   JOURNAL = {Math. Program.},
  FJOURNAL = {Mathematical Programming},
    VOLUME = {151},
      YEAR = {2015},
    NUMBER = {1, Ser. B},
     PAGES = {3--34},
      ISSN = {0025-5610},
   MRCLASS = {90C25 (49M20 65K05)},
  MRNUMBER = {3347548},
MRREVIEWER = {Marco Gaviano},
       DOI = {10.1007/s10107-015-0892-3},
       URL = {https://doi.org/10.1007/s10107-015-0892-3},
}

@article{CDPR2003,
AUTHOR = {Canutescu, Adrian A. and Dunbrack, Roland L.  },
TITLE =  {Cyclic coordinate descent: {A} robotics algorithm for protein loop closure},
JOURNAL = {Protein Sci.},
VOLUME = {12},
PAGES = {963--972},
YEAR = {2003},
NUMBER = {5},
}

@conference{Liu-2009-10238,
author = {Han Liu and Mark Palatucci and Jian Zhang},
title = {Blockwise Coordinate Descent Procedures for the Multi-task Lasso, with Applications to Neural Semantic Basis Discovery},
booktitle = {Proceedings of (ICML) International Conference on Machine Learning},
year = {2009},
month = {June},
pages = {649 - 656},
}

@article{06ae6e67-7fa9-32aa-97a2-a282f021e2c2,
 ISSN = {10618600},
 URL = {http://www.jstor.org/stable/1390659},
 author = {Sylvain Sardy and Andrew G. Bruce and Paul Tseng},
 journal = {J. Comput.  Graph. Stat.},
 number = {2},
 pages = {361--379},
 publisher = {[American Statistical Association, Taylor & Francis, Ltd., Institute of Mathematical Statistics, Interface Foundation of America]},
 title = {Block Coordinate Relaxation Methods for Nonparametric Wavelet Denoising},
 urldate = {2025-03-19},
 volume = {9},
 year = {2000}
}

@article {MR3179953,
    AUTHOR = {Richt\'{a}rik, Peter and Tak\'{a}\v{c}, Martin},
     TITLE = {Iteration complexity of randomized block-coordinate descent
              methods for minimizing a composite function},
   JOURNAL = {Math. Program.},
  FJOURNAL = {Mathematical Programming},
    VOLUME = {144},
      YEAR = {2014},
    NUMBER = {1-2, Ser. A},
     PAGES = {1--38},
      ISSN = {0025-5610},
   MRCLASS = {90C06 (65K10 65Y20 90C25 90C60)},
  MRNUMBER = {3179953},
       DOI = {10.1007/s10107-012-0614-z},
       URL = {https://doi.org/10.1007/s10107-012-0614-z},
}

@inproceedings{10.5555/3104482.3104523,
author = {Bradley, Joseph K. and Kyrola, Aapo and Bickson, Danny and Guestrin, Carlos},
title = {Parallel coordinate descent for {L}1-regularized loss minimization},
year = {2011},
isbn = {9781450306195},
publisher = {Omnipress},
address = {Madison, WI, USA},
booktitle = {Proceedings of the 28th International Conference on International Conference on Machine Learning},
pages = {321–328},
numpages = {8},
location = {Bellevue, Washington, USA},
series = {ICML'11}
}

@article {MR3459207,
    AUTHOR = {Richt\'{a}rik, Peter and Tak\'{a}\v{c}, Martin},
     TITLE = {Parallel coordinate descent methods for big data optimization},
   JOURNAL = {Math. Program.},
  FJOURNAL = {Mathematical Programming},
    VOLUME = {156},
      YEAR = {2016},
    NUMBER = {1-2, Ser. A},
     PAGES = {433--484},
      ISSN = {0025-5610},
   MRCLASS = {90C06 (49M20 49M27 65K05 68W10 68W20 68W40 90C25)},
  MRNUMBER = {3459207},
       DOI = {10.1007/s10107-015-0901-6},
       URL = {https://doi.org/10.1007/s10107-015-0901-6},
}

@article {MR3513271,
    AUTHOR = {Tappenden, Rachael and Richt\'{a}rik, Peter and Gondzio, Jacek},
     TITLE = {Inexact coordinate descent: complexity and preconditioning},
   JOURNAL = {J. Optim. Theory Appl.},
  FJOURNAL = {Journal of Optimization Theory and Applications},
    VOLUME = {170},
      YEAR = {2016},
    NUMBER = {1},
     PAGES = {144--176},
      ISSN = {0022-3239},
   MRCLASS = {65F08 (65F10 65K10 65Y20 68Q25 90C25)},
  MRNUMBER = {3513271},
MRREVIEWER = {Arvind Krishna Saibaba},
       DOI = {10.1007/s10957-016-0867-4},
       URL = {https://doi.org/10.1007/s10957-016-0867-4},
}

@article {MR3369495,
    AUTHOR = {Lu, Zhaosong and Xiao, Lin},
     TITLE = {On the complexity analysis of randomized block-coordinate
              descent methods},
   JOURNAL = {Math. Program.},
  FJOURNAL = {Mathematical Programming},
    VOLUME = {152},
      YEAR = {2015},
    NUMBER = {1-2, Ser. A},
     PAGES = {615--642},
      ISSN = {0025-5610},
   MRCLASS = {65K05 (65Y20 90C05 90C06 90C25)},
  MRNUMBER = {3369495},
MRREVIEWER = {Nicolas Gillis},
       DOI = {10.1007/s10107-014-0800-2},
       URL = {https://doi.org/10.1007/s10107-014-0800-2},
}

@article {MR3517098,
    AUTHOR = {Richt\'{a}rik, Peter and Tak\'{a}\v{c}, Martin},
     TITLE = {Distributed coordinate descent method for learning with big
              data},
   JOURNAL = {J. Mach. Learn. Res.},
  FJOURNAL = {Journal of Machine Learning Research (JMLR)},
    VOLUME = {17},
      YEAR = {2016},
     PAGES = {Paper No. 75, 25},
      ISSN = {1532-4435},
   MRCLASS = {68T05 (62H30 68W15)},
  MRNUMBER = {3517098},
}

@article {MR321541,
    AUTHOR = {Powell, M. J. D.},
     TITLE = {On search directions for minimization algorithms},
   JOURNAL = {Math. Program.},
  FJOURNAL = {Mathematical Programming},
    VOLUME = {4},
      YEAR = {1973},
     PAGES = {193--201},
      ISSN = {0025-5610},
   MRCLASS = {90C30},
  MRNUMBER = {321541},
MRREVIEWER = {K. Zimmermann},
       DOI = {10.1007/BF01584660},
       URL = {https://doi.org/10.1007/BF01584660},
}

@book {MR1955649,
	AUTHOR = {Facchinei, F. and Pang, J.-S.},
	TITLE = {Finite-Dimensional Variational Inequalities and
	Complementarity Problems, {V}olume {II}},
	PUBLISHER = {Springer-Verlag, New York},
	YEAR = {2003},
	PAGES = {i-xxxiv, 625--1234 and II1--II57},
	ISBN = {0-387-95581-X},
	MRCLASS = {90-02 (15A39 49J40 65K10 90C33)},
	MRNUMBER = {1955649},
	MRREVIEWER = {Daniel Ralph},
}

@article {MR3365070,
    AUTHOR = {Patrascu, Andrei and Necoara, Ion},
     TITLE = {Random coordinate descent methods for {$\ell_0$} regularized
              convex optimization},
   JOURNAL = {IEEE Trans. Automat. Control},
  FJOURNAL = {Institute of Electrical and Electronics Engineers.
              Transactions on Automatic Control},
    VOLUME = {60},
      YEAR = {2015},
    NUMBER = {7},
     PAGES = {1811--1824},
      ISSN = {0018-9286,1558-2523},
   MRCLASS = {90C26 (49M37 90C46 90C59 94C99)},
  MRNUMBER = {3365070},
MRREVIEWER = {Driss\ Mentagui},
       DOI = {10.1109/TAC.2015.2390551},
       URL = {https://doi.org/10.1109/TAC.2015.2390551},
}

@INPROCEEDINGS{9414191,
  author={Gao, Tianxiang and Lu, Songtao and Liu, Jia and Chu, Chris},
  booktitle={ICASSP 2021 - 2021 IEEE International Conference on Acoustics, Speech and Signal Processing (ICASSP)}, 
  title={On the Convergence of Randomized Bregman Coordinate Descent for Non-{L}ipschitz Composite Problems}, 
  year={2021},
  volume={},
  number={},
  pages={5549-5553},
  doi={10.1109/ICASSP39728.2021.9414191}}

@article {MR2421270,
    AUTHOR = {Attouch, Hedy and Bolte, J\'{e}r\^{o}me},
     TITLE = {On the convergence of the proximal algorithm for nonsmooth
              functions involving analytic features},
   JOURNAL = {Math. Program.},
  FJOURNAL = {Mathematical Programming. A Publication of the Mathematical
              Programming Society},
    VOLUME = {116},
      YEAR = {2009},
    NUMBER = {1-2, Ser. B},
     PAGES = {5--16},
      ISSN = {0025-5610},
   MRCLASS = {90C30 (65K05)},
  MRNUMBER = {2421270},
MRREVIEWER = {W. W. Breckner},
       DOI = {10.1007/s10107-007-0133-5},
       URL = {https://doi.org/10.1007/s10107-007-0133-5},
}

@article {MR3010421,
    AUTHOR = {Attouch, Hedy and Bolte, J\'{e}r\^{o}me and Svaiter, Benar Fux},
     TITLE = {Convergence of descent methods for semi-algebraic and tame
              problems: proximal algorithms, forward-backward splitting, and
              regularized {G}auss-{S}eidel methods},
   JOURNAL = {Math. Program.},
  FJOURNAL = {Mathematical Programming},
    VOLUME = {137},
      YEAR = {2013},
    NUMBER = {1-2, Ser. A},
     PAGES = {91--129},
      ISSN = {0025-5610},
   MRCLASS = {49J53 (47J25 47J30 49M15 65K10 90C30)},
  MRNUMBER = {3010421},
MRREVIEWER = {C. H. Jeffrey Pang},
       DOI = {10.1007/s10107-011-0484-9},
       URL = {https://doi.org/10.1007/s10107-011-0484-9},
}

@article {MR3232623,
    AUTHOR = {Bolte, J\'{e}r\^{o}me and Sabach, Shoham and Teboulle, Marc},
     TITLE = {Proximal alternating linearized minimization for nonconvex and
              nonsmooth problems},
   JOURNAL = {Math. Program.},
  FJOURNAL = {Mathematical Programming},
    VOLUME = {146},
      YEAR = {2014},
    NUMBER = {1-2, Ser. A},
     PAGES = {459--494},
      ISSN = {0025-5610},
   MRCLASS = {90C26 (90C56)},
  MRNUMBER = {3232623},
MRREVIEWER = {Yiran He},
       DOI = {10.1007/s10107-013-0701-9},
       URL = {https://doi.org/10.1007/s10107-013-0701-9},
}

@article {MR2674728,
    AUTHOR = {Attouch, H\'{e}dy and Bolte, J\'{e}r\^{o}me and Redont, Patrick and
              Soubeyran, Antoine},
     TITLE = {Proximal alternating minimization and projection methods for
              nonconvex problems: an approach based on the
              {K}urdyka-{{\L}ojasiewicz} inequality},
   JOURNAL = {Math. Oper. Res.},
  FJOURNAL = {Mathematics of Operations Research},
    VOLUME = {35},
      YEAR = {2010},
    NUMBER = {2},
     PAGES = {438--457},
      ISSN = {0364-765X},
   MRCLASS = {90C26 (49J52 65K10)},
  MRNUMBER = {2674728},
       DOI = {10.1287/moor.1100.0449},
       URL = {https://doi.org/10.1287/moor.1100.0449},
}

@article {MR3218822,
    AUTHOR = {Ochs, Peter and Chen, Yunjin and Brox, Thomas and Pock,
              Thomas},
     TITLE = {i{P}iano: inertial proximal algorithm for nonconvex
              optimization},
   JOURNAL = {SIAM J. Imaging Sci.},
  FJOURNAL = {SIAM Journal on Imaging Sciences},
    VOLUME = {7},
      YEAR = {2014},
    NUMBER = {2},
     PAGES = {1388--1419},
   MRCLASS = {94A08},
  MRNUMBER = {3218822},
       DOI = {10.1137/130942954},
       URL = {https://doi.org/10.1137/130942954},
}

@article {MR3327417,
    AUTHOR = {Bredies, Kristian and Lorenz, Dirk A. and Reiterer, Stefan},
     TITLE = {Minimization of non-smooth, non-convex functionals by
              iterative thresholding},
   JOURNAL = {J. Optim. Theory Appl.},
  FJOURNAL = {Journal of Optimization Theory and Applications},
    VOLUME = {165},
      YEAR = {2015},
    NUMBER = {1},
     PAGES = {78--112},
      ISSN = {0022-3239},
   MRCLASS = {49M05 (65K10 90C26 90C48)},
  MRNUMBER = {3327417},
MRREVIEWER = {Jan Mach},
       DOI = {10.1007/s10957-014-0614-7},
       URL = {https://doi.org/10.1007/s10957-014-0614-7},
}

@article {MR4690351,
    AUTHOR = {Yang, Lei},
     TITLE = {Proximal gradient method with extrapolation and line search
              for a class of non-convex and non-smooth problems},
   JOURNAL = {J. Optim. Theory Appl.},
  FJOURNAL = {Journal of Optimization Theory and Applications},
    VOLUME = {200},
      YEAR = {2024},
    NUMBER = {1},
     PAGES = {68--103},
      ISSN = {0022-3239},
   MRCLASS = {90C26 (65K10)},
  MRNUMBER = {4690351},
       DOI = {10.1007/s10957-023-02348-4},
       URL = {https://doi.org/10.1007/s10957-023-02348-4},
}

@INPROCEEDINGS{8263933,
  author={Stella, Lorenzo and Themelis, Andreas and Sopasakis, Pantelis and Patrinos, Panagiotis},
  booktitle={2017 IEEE 56th Annual Conference on Decision and Control (CDC)}, 
  title={A simple and efficient algorithm for nonlinear model predictive control}, 
  year={2017},
  volume={},
  number={},
  pages={1939-1944},
  doi={10.1109/CDC.2017.8263933}}

@article{MR4457955,
    AUTHOR = {De Marchi, Alberto and Themelis, Andreas},
     TITLE = {Proximal gradient algorithms under local {L}ipschitz gradient
              continuity: a convergence and robustness analysis of {PANOC}},
   JOURNAL = {J. Optim. Theory Appl.},
  FJOURNAL = {Journal of Optimization Theory and Applications},
    VOLUME = {194},
      YEAR = {2022},
    NUMBER = {3},
     PAGES = {771--794},
      ISSN = {0022-3239},
   MRCLASS = {90C26 (49J52 65K05)},
  MRNUMBER = {4457955},
MRREVIEWER = {Najmeh Hoseini Monjezi},
       DOI = {10.1007/s10957-022-02048-5},
       URL = {https://doi.org/10.1007/s10957-022-02048-5},
}

@article {MR4665669,
    AUTHOR = {Jia, Xiaoxi and Kanzow, Christian and Mehlitz, Patrick},
     TITLE = {Convergence analysis of the proximal gradient method in the
              presence of the {K}urdyka-{{\L}ojasiewicz} property without global
              {L}ipschitz assumptions},
   JOURNAL = {SIAM J. Optim.},
  FJOURNAL = {SIAM Journal on Optimization},
    VOLUME = {33},
      YEAR = {2023},
    NUMBER = {4},
     PAGES = {3038--3056},
      ISSN = {1052-6234},
   MRCLASS = {90C30 (49J52)},
  MRNUMBER = {4665669},
       DOI = {10.1137/23M1548293},
       URL = {https://doi.org/10.1137/23M1548293},
}

@article {MR4927930,
    AUTHOR = {Kanzow, Christian and Lehmann, Leo},
     TITLE = {Convergence of nonmonotone proximal gradient methods under the
              {K}urdyka--{{\L}ojasiewicz} property without a global {L}ipschitz
              assumption},
   JOURNAL = {J. Optim. Theory Appl.},
  FJOURNAL = {Journal of Optimization Theory and Applications},
    VOLUME = {207},
      YEAR = {2025},
    NUMBER = {1},
     PAGES = {Paper No. 4, 26},
      ISSN = {0022-3239},
   MRCLASS = {90C30 (49J52 65K05)},
  MRNUMBER = {4927930},
       DOI = {10.1007/s10957-025-02762-w},
       URL = {https://doi.org/10.1007/s10957-025-02762-w},
}

@inproceedings{LeeSeung2000,
 author = {Lee, Daniel and Seung, H. Sebastian},
 booktitle = {Advances in Neural Information Processing Systems},
 editor = {T. Leen and T. Dietterich and V. Tresp},
 pages = {},
 publisher = {MIT Press},
 title = {Algorithms for Non-negative Matrix Factorization},
 url = {https://proceedings.neurips.cc/paper_files/paper/2000/file/f9d1152547c0bde01830b7e8bd60024c-Paper.pdf},
 volume = {13},
 year = {2000}
}

@inproceedings{DonohoStodden2003,
 author = {Donoho, David and Stodden, Victoria},
 booktitle = {Advances in Neural Information Processing Systems},
 editor = {S. Thrun and L. Saul and B. Sch\"{o}lkopf},
 pages = {},
 publisher = {MIT Press},
 title = {When Does Non-Negative Matrix Factorization Give a Correct Decomposition into Parts?},
 url = {https://proceedings.neurips.cc/paper_files/paper/2003/file/1843e35d41ccf6e63273495ba42df3c1-Paper.pdf},
 volume = {16},
 year = {2003}
}

@article {MR3572363,
    AUTHOR = {Pock, Thomas and Sabach, Shoham},
     TITLE = {Inertial proximal alternating linearized minimization
              (i{PALM}) for nonconvex and nonsmooth problems},
   JOURNAL = {SIAM J. Imaging Sci.},
  FJOURNAL = {SIAM Journal on Imaging Sciences},
    VOLUME = {9},
      YEAR = {2016},
    NUMBER = {4},
     PAGES = {1756--1787},
      ISSN = {1936-4954},
   MRCLASS = {94A08 (49M37 65K10 90C26 90C30)},
  MRNUMBER = {3572363},
MRREVIEWER = {A.\ G.\ Kartsatos},
       DOI = {10.1137/16M1064064},
       URL = {https://doi.org/10.1137/16M1064064},
}

@article {MR4354998,
    AUTHOR = {Driggs, Derek and Tang, Junqi and Liang, Jingwei and Davies,
              Mike and Sch\"{o}nlieb, Carola-Bibiane},
     TITLE = {A stochastic proximal alternating minimization for nonsmooth
              and nonconvex optimization},
   JOURNAL = {SIAM J. Imaging Sci.},
  FJOURNAL = {SIAM Journal on Imaging Sciences},
    VOLUME = {14},
      YEAR = {2021},
    NUMBER = {4},
     PAGES = {1932--1970},
      ISSN = {1936-4954},
   MRCLASS = {90C26 (90C06 90C15 90C56)},
  MRNUMBER = {4354998},
MRREVIEWER = {Najmeh\ Hoseini Monjezi},
       DOI = {10.1137/20M1387213},
       URL = {https://doi.org/10.1137/20M1387213},
}

@article {MR3353930,
    AUTHOR = {Kuang, Da and Yun, Sangwoon and Park, Haesun},
     TITLE = {Sym{NMF}: nonnegative low-rank approximation of a similarity
              matrix for graph clustering},
   JOURNAL = {J. Global Optim.},
  FJOURNAL = {Journal of Global Optimization. An International Journal
              Dealing with Theoretical and Computational Aspects of Seeking
              Global Optima and Their Applications in Science, Management
              and Engineering},
    VOLUME = {62},
      YEAR = {2015},
    NUMBER = {3},
     PAGES = {545--574},
      ISSN = {0925-5001,1573-2916},
   MRCLASS = {62H30 (05C82 15A23 15B48 62B10 65F99 90C90)},
  MRNUMBER = {3353930},
       DOI = {10.1007/s10898-014-0247-2},
       URL = {https://doi.org/10.1007/s10898-014-0247-2},
}

@inproceedings{DingSymNMF,
author = {Chris Ding and Xiaofeng He and Horst D. Simon},
title = {On the Equivalence of Nonnegative Matrix Factorization and Spectral Clustering},
booktitle = {Proceedings of the 2005 SIAM International Conference on Data Mining (SDM)},
chapter = {},
pages = {606-610},
doi = {10.1137/1.9781611972757.70},
URL = {https://epubs.siam.org/doi/abs/10.1137/1.9781611972757.70},
eprint = {https://epubs.siam.org/doi/pdf/10.1137/1.9781611972757.70}
}

@article {MR2785131,
    AUTHOR = {Wang, Fei and Li, Tao and Wang, Xin and Zhu, Shenghuo and
              Ding, Chris},
     TITLE = {Community discovery using nonnegative matrix factorization},
   JOURNAL = {Data Min. Knowl. Discov.},
  FJOURNAL = {Data Mining and Knowledge Discovery},
    VOLUME = {22},
      YEAR = {2011},
    NUMBER = {3},
     PAGES = {493--521},
      ISSN = {1384-5810,1573-756X},
   MRCLASS = {68M10 (05C82 90B10)},
  MRNUMBER = {2785131},
       DOI = {10.1007/s10618-010-0181-y},
       URL = {https://doi.org/10.1007/s10618-010-0181-y},
}

@article{scikit,
author = {Pedregosa, Fabian and Varoquaux, Ga\"{e}l and Gramfort, Alexandre and Michel, Vincent and Thirion, Bertrand and Grisel, Olivier and Blondel, Mathieu and Prettenhofer, Peter and Weiss, Ron and Dubourg, Vincent and Vanderplas, Jake and Passos, Alexandre and Cournapeau, David and Brucher, Matthieu and Perrot, Matthieu and Duchesnay, \'{E}douard},
title = {Scikit-learn: {M}achine {L}earning in {P}ython},
year = {2011},
issue_date = {2/1/2011},
publisher = {JMLR.org},
volume = {12},
number = {null},
issn = {1532-4435},
journal = {J. Mach. Learn. Res.},
month = nov,
pages = {2825--2830},
numpages = {6}
}

@article {MR2409803,
    AUTHOR = {von Luxburg, Ulrike},
     TITLE = {A tutorial on spectral clustering},
   JOURNAL = {Stat. Comput.},
  FJOURNAL = {Statistics and Computing},
    VOLUME = {17},
      YEAR = {2007},
    NUMBER = {4},
     PAGES = {395--416},
      ISSN = {0960-3174,1573-1375},
   MRCLASS = {Expansion},
  MRNUMBER = {2409803},
       DOI = {10.1007/s11222-007-9033-z},
       URL = {https://doi.org/10.1007/s11222-007-9033-z},
}

\end{document}